# ATLAS MODELS OF EQUITY MARKETS

By Adrian D. Banner, Robert Fernholz and Ioannis Karatzas

*INTECH, INTECH and Columbia University*

Atlas-type models are constant-parameter models of uncorrelated stocks for equity markets with a stable capital distribution, in which the growth rates and variances depend on rank. The simplest such model assigns the same, constant variance to all stocks; zero rate of growth to all stocks but the smallest; and positive growth rate to the smallest, the *Atlas* stock. In this paper we study the basic properties of this class of models, as well as the behavior of various portfolios in their midst. Of particular interest are portfolios that do not contain the Atlas stock.

**1. Introduction.** Size is one of the most important descriptive characteristics of assets: one can understand a great deal about an equity market by observing, and making sense of, the continual ebb and flow of small-, medium- and large-capitalization stocks in its midst. Thus it is important to have models which describe (if not explain) this flow, and which exhibit stability properties for the resulting distribution of capital that are in agreement with actual observation. This paper studies models of this type and analyzes portfolio performance in their context.

The simplest such model is the *Atlas model* for equity markets, introduced in Example 5.3.3 of [2]. This is a constant-coefficient model for the values (capitalizations) of stocks represented by their relative rank and driven by independent Brownian motions. It assigns the same, constant volatility to all stocks; zero growth rate to all stocks but the smallest; and positive growth rate to the smallest stock. Because it is responsible for all the growth (or support) in the market, this smallest stock is then called the *Atlas stock*.

Somewhat more precisely: with $g > 0$, $\sigma > 0$ given constants, with independent Brownian motions $W_1(\cdot), \ldots, W_n(\cdot)$, and with $X_i(t)$ representing the capitalization at time $t$ of the stock with index (name) $i$, the Atlas model

Received November 2004.

---

postulates the dynamics

$$(1.1) \qquad d(\log X_i(t)) = \gamma_i(t)\,dt + \sigma_i(t)\,dW_i(t), \qquad i = 1, \ldots, n,$$

where the growth rates and volatilities are specified by

$$(1.2) \qquad \gamma_i(t) = ng \cdot \mathbb{1}_{\{X_i(t) = X_{p_t(n)}(t)\}}, \qquad \sigma_i(t) = \sigma.$$

We are using here the "reverse order-statistics" notation

$$(1.3) \quad \begin{aligned} \max_{1 \le i \le n} X_i(t) =: X_{(1)}(t) \ge X_{(2)}(t) \ge \\ \cdots \ge X_{(n-1)}(t) \ge X_{(n)}(t) := \min_{1 \le i \le n} X_i(t) \end{aligned}$$

for the capitalizations of stocks ranked in descending order, from largest to smallest; we consider also the random permutation $(p_t(1), \ldots, p_t(n))$ of $(1, \ldots, n)$, for which

$$(1.4) \qquad X_{p_t(k)}(t) = X_{(k)}(t), \qquad p_t(k) < p_t(k+1) \text{ if } X_{(k)}(t) = X_{(k+1)}(t)$$

hold with $k = 1, \ldots, n$. Roughly speaking, this means that $p_t(k)$ is the name (index) of the stock with the $k$th-largest relative capitalization at time $t$, and that ties are resolved by resorting to the lowest index.

More generally, suppose we are given real numbers $\gamma$, $g_1, \ldots, g_n$, $\sigma_1 > 0, \ldots, \sigma_n > 0$ such that

$$(1.5) \qquad \begin{aligned} g_1 < 0, \qquad g_1 + g_2 < 0, \quad \ldots, \quad g_1 + \cdots + g_{n-1} < 0, \\ g_1 + \cdots + g_n = 0. \end{aligned}$$

Corresponding to these parameters, the general model considered in this paper postulates the dynamics of (1.1) for the stock capitalizations $X_1(t), \ldots, X_n(t)$, but now with growth rates and volatilities given by

$$(1.6) \quad \gamma_i(t) = \gamma + \sum_{k=1}^{n} g_k \mathbb{1}_{\{X_i(t) = X_{p_t(k)}(t)\}}, \qquad \sigma_i(t) = \sum_{k=1}^{n} \sigma_k \mathbb{1}_{\{X_i(t) = X_{p_t(k)}(t)\}}$$

in place of (1.2). In other words, this more general model specifies $\gamma + g_k$ as the growth rate, and $\sigma_k$ as the volatility, for the stock with rank $k$ at any given time. We shall refer to the model of (1.1), (1.6) as the *first-order model*.

Clearly

$$(1.7) \qquad \gamma = g > 0, \qquad g_k = -g \quad \text{for } k = 1, \ldots, n-1 \quad \text{and} \quad g_n = (n-1)g,$$

in the case of the Atlas model of (1.1), (1.2). We shall call *generalized Atlas model* a model of the type (1.1), (1.6) with parameters that satisfy (1.7), though with possibly different volatilities.

All these models have strictly nondegenerate volatility structures and bounded drift coefficients, so they admit a (unique) equivalent martingale measure on any given time-horizon; thus *there are no relative arbitrage opportunities for such models*, of the type encountered in [3, 4].

The first question that arises for the first-order model of (1.1), (1.6) is the rigorous formulation of, and the study of existence/uniqueness of solution to, the resulting system of stochastic differential equations. This task we undertake in Section 2, whereas in Section 3 and the Appendix we study the behavior of the resulting "ranked capitalization" (reverse-order-statistics) processes of (1.3). Section 4 deals with ergodic properties of these processes. Portfolios in the context of the model (1.1), (1.6) are introduced in Section 5, where we also study the growth rates of a few relatively easy-to-implement investment rules. Some detailed comparisons of long-term-growth performance are carried out in Section 6. We conclude with considerations of diversity in Section 7, where some elementary computations show that models of this sort capture very well the intuitive notion that "no stock can be allowed to dominate the entire market with anything but extremely low probability"—despite the fact that such models fail to be diverse in a strict, almost sure sense.

**2. The model.** Let us start by constructing a diffusion process corresponding to the stochastic equation of (1.1), (1.6). We consider a collection $\{\mathcal{Q}_k^{(i)}\}_{1 \le i, k \le n}$ of polyhedral domains in $\mathbb{R}^n$ with the following property:

(2.1) $y = (y_1, \ldots, y_n) \in \mathcal{Q}_k^{(i)}$ means that $y_i$ is ranked $k$th among $y_1, \ldots, y_n$.

We resolve ties by resorting to the lowest index; for instance, we set

$$\mathcal{Q}_1^{(1)} = \{y \in \mathbb{R}^n | y_1 \ge y_j \ \forall j = 2, \ldots, n\},$$

$$\mathcal{Q}_n^{(1)} = \{y \in \mathbb{R}^n | y_1 < y_j \ \forall j = 2, \ldots, n\}$$

and

$$\mathcal{Q}_{k+1}^{(1)} = \left\{y \in \mathbb{R}^n \,\middle|\, y_1 < \min_{1 \le r \le k} y_{j_r} \text{ for some } j_1, \ldots, j_k, \right.$$
$$\left. \text{and } y_1 \ge y_\ell \ \forall \ell \notin \{j_1, \ldots, j_k\} \right\}$$

for $k = 1, \ldots, n - 2$. Clearly, the collection $\{\mathcal{Q}_k^{(i)}\}_{1 \le i \le n}$ is a partition of $\mathbb{R}^n$ for each fixed $k$; and

(2.2) $\{\mathcal{Q}_k^{(i)}\}_{1 \le k \le n}$ is a partition of $\mathbb{R}^n$ for each fixed $i$.

Consider now real constants $\gamma, \sigma_1 > 0, \ldots, \sigma_n > 0$ and $g_1, \ldots, g_n$ satisfying the conditions of (1.5). We shall look at the system of stochastic differential

equations

$$
\begin{aligned}
(2.3) \quad dY_i(t) = {} & \left( \sum_{k=1}^{n} g_k \mathbb{1}_{\mathcal{Q}_k^{(i)}}(Y(t)) + \gamma \right) dt \\
& + \sum_{k=1}^{n} \sigma_k \mathbb{1}_{\mathcal{Q}_k^{(i)}}(Y(t)) \cdot dW_i(t), \qquad Y_i(0) = \bar{y}_i,
\end{aligned}
$$

for the $n$-dimensional process $Y(\cdot) = (Y_1(\cdot), \ldots, Y_n(\cdot))$, with given initial condition $\bar{y} = (\bar{y}_1, \ldots, \bar{y}_n)$. In other words: as long as $Y(\cdot)$ is in the polyhedron $\mathcal{Q}_k^{(i)}$, equation (2.3) postulates that the $i$th-coordinate process $Y_i(\cdot)$ evolve like a Brownian motion with drift $g_k + \gamma$ and variance $\sigma_k^2$, for each $i = 1, \ldots, n$.

The theory of Bass and Pardoux [1] establishes that the system of stochastic differential equations (2.3) has a weak solution, which is unique in the sense of the probability law. Once this solution has been constructed, we can look at the vector of processes $X(\cdot) = (X_1(\cdot), \ldots, X_n(\cdot))$ defined by

$$
(2.4) \qquad X_i(t) := e^{Y_i(t)}, \qquad i = 1, \ldots, n,
$$

as the rigorous interpretation of the first-order model of (1.1), (1.6). With this interpretation $Y_i(t)$ represents the log-capitalization of the $i$th company at time $t$.

Now let us observe from (2.3), the remark preceding (2.2), and (1.5), that we have

$$
\begin{aligned}
(2.5) \quad \sum_{i=1}^{n} Y_i(t) = {} & Y + n\gamma t + \sum_{k=1}^{n} \sigma_k B_k(t) \\
& \text{where } B_k(t) := \sum_{i=1}^{n} \int_0^t \mathbb{1}_{\mathcal{Q}_k^{(i)}}(Y(s)) \, dW_i(s)
\end{aligned}
$$

and $Y = \sum_{i=1}^{n} Y_i(0)$. The resulting processes $B_1(\cdot), \ldots, B_n(\cdot)$ are continuous local martingales with quadratic (cross-) variations $\langle B_k, B_\ell \rangle(t)$ equal to

$$
\sum_{i=1}^{n} \sum_{j=1}^{n} \int_0^t \mathbb{1}_{\mathcal{Q}_k^{(i)}}(Y(s)) \mathbb{1}_{\mathcal{Q}_\ell^{(j)}}(Y(s)) \, d\langle W_i, W_j \rangle(s)
$$

$$
= \sum_{i=1}^{n} \int_0^t \mathbb{1}_{\mathcal{Q}_k^{(i)}}(Y(s)) \mathbb{1}_{\mathcal{Q}_\ell^{(i)}}(Y(s)) \, ds = \delta_{k\ell} t.
$$

Lévy's characterization (e.g., Theorem 3.3.16 in [8]) identifies the processes $B_1(\cdot), \ldots, B_n(\cdot)$ as independent standard Brownian motions; then the strong law of large numbers in conjunction with (2.5) gives

$$
(2.6) \qquad \lim_{T \to \infty} \frac{1}{T} \sum_{i=1}^{n} Y_i(T) = n\gamma \qquad \text{a.s.}
$$

In the Appendix we shall strengthen this result, and show that in fact

$$(2.7) \qquad \lim_{T \to \infty} \frac{1}{T} \log X_i(T) = \lim_{T \to \infty} \frac{Y_i(T)}{T} = \gamma$$

holds a.s. for every $i = 1, \dots, n$.

REMARK 2.1 (Coherence).   Denoting by $X(t) := X_1(t) + \cdots + X_n(t)$ the total market capitalization and by

$$(2.8) \qquad \mu_i(t) := \frac{X_i(t)}{X(t)}, \qquad i = 1, \dots, n,$$

the relative capitalizations of the individual companies, we see from (2.7) that

$$\lim_{T \to \infty} \frac{1}{T} \log X(T) = \max_{1 \le i \le n} \left( \lim_{T \to \infty} \frac{1}{T} \log X_i(T) \right) = \gamma,$$

and thus also

$$\lim_{T \to \infty} \frac{1}{T} \log \mu_i(T) = 0 \qquad \forall\, i = 1, \dots, n,$$

holds a.s. In the terminology of Fernholz ([2], page 26), the model of (2.3), (2.4)—or equivalently, that of (1.1), (1.6)—is *coherent*.

REMARK 2.2 (Taking turns as Atlas).   From (2.3), (2.7) and the strong law of large numbers for Brownian motion, it follows that

$$\lim_{T \to \infty} \sum_{k=1}^{n} g_k \left( \frac{1}{T} \int_0^T \mathbb{1}_{\mathcal{Q}_k^{(i)}}(Y(t))\, dt \right) = 0 \qquad \text{holds a.s.}$$

for every $i = 1, \dots, n$. Now suppose the parameters of the model satisfy the conditions $\gamma = g > 0$, $g_1 = \cdots = g_{n-1} = -g$, $g_n = (n-1)g$ of (1.7) for a generalized Atlas model; then

$$\sum_{k=1}^{n} g_k \left( \frac{1}{T} \int_0^T \mathbb{1}_{\mathcal{Q}_k^{(i)}}(Y(t))\, dt \right) = g \left( \frac{n}{T} \int_0^T \mathbb{1}_{\mathcal{Q}_n^{(i)}}(Y(t))\, dt - 1 \right),$$

and for every $i = 1, \dots, n$ we obtain

$$(2.9) \qquad \lim_{T \to \infty} \frac{1}{T} \int_0^T \mathbb{1}_{\mathcal{Q}_n^{(i)}}(Y(t))\, dt = \frac{1}{n} \qquad \text{a.s.}$$

In other words: "each stock acts as Atlas roughly $(1/n)$th of the time." It is then natural to conjecture that we should have

$$(2.10) \qquad \lim_{T \to \infty} \frac{1}{T} \int_0^T \mathbb{1}_{\mathcal{Q}_k^{(i)}}(Y(t))\, dt = \frac{1}{n} \qquad \text{a.s.}$$

*for every* $k = 1, \dots, n$, that is, not just for $k = n$ as in (2.9).

As it turns out, this property holds for the general first-order model; in particular, each stock spends asymptotically the same amount of time in every rank.

PROPOSITION 2.3. *The solution of the system* (2.3) *of stochastic differential equations satisfies the ergodic relation* (2.10) *for every* $k = 1, \ldots, n$ *and every* $i = 1, \ldots, n$.

PROOF. Let $\Sigma_n$ denote the symmetric group of permutations of $\{1, \ldots, n\}$. For each $p \in \Sigma_n$, let $\mathcal{R}_p := \bigcap_{k=1}^{n} \mathcal{Q}_k^{(p(k))}$; the set $\mathcal{R}_p$ consists of all points $y = (y_1, \ldots, y_n) \in \mathbb{R}^n$ such that $y_{p(k)}$ is ranked $k$th among $y_1, \ldots, y_n$ for all $k = 1, \ldots, n$ (with ties once again resolved in favor of the lowest index). In particular,

$$(2.11) \qquad y_{p(1)} \geq y_{p(2)} \geq \cdots \geq y_{p(n)} \qquad \text{if } y \in \mathcal{R}_p.$$

Clearly $\{\mathcal{R}_p\}_{p \in \Sigma_n}$ is a partition of $\mathbb{R}^n$. Let $\mathcal{R}_p^\circ$ denote the interior of the polyhedron $\mathcal{R}_p$, and set $\mathcal{R}^\circ := \bigcup_{p \in \Sigma_n} \mathcal{R}_p^\circ$. The exceptional set $\mathcal{R}^e := \mathbb{R}^n \backslash \mathcal{R}^\circ$ can also be described as $\{(y_1, \ldots, y_n) \in \mathbb{R}^n | y_i = y_j \text{ for some } i \neq j\}$. Furthermore, any $p \in \Sigma_n$ acts as a linear transformation of $\mathbb{R}^n$ via $p(y_1, \ldots, y_n) := (y_{p^{-1}(1)}, \ldots, y_{p^{-1}(n)})$; under this action, we have

$$(2.12) \quad y \in \mathcal{R}_p^\circ \quad \Longleftrightarrow \quad p^{-1}y \in \mathcal{R}_1^\circ := \{(y_1, \ldots, y_n) \in \mathbb{R}^n | y_1 > y_2 > \cdots > y_n\}.$$

Define

$$G(y) := \sum_{p \in \Sigma_n} \mathbb{1}_{\mathcal{R}_p}(y)(g_{p^{-1}(1)}, \ldots, g_{p^{-1}(n)})^t,$$

$$S(y) := \sum_{p \in \Sigma_n} \mathbb{1}_{\mathcal{R}_p}(y) \operatorname{diag}(\sigma_{p^{-1}(1)}, \ldots, \sigma_{p^{-1}(n)}).$$

Set $v := (1, 1, \ldots, 1)^t$ and note that $y \in \mathcal{R}_p \Longleftrightarrow y + \alpha v \in \mathcal{R}_p$ for all $\alpha \in \mathbb{R}$, as ranks of coordinates are preserved by adding scalar multiples of $v$. It follows that

$$(2.13) \quad G(y + \alpha v) = G(y) \quad \text{and} \quad S(y + \alpha v) = S(y) \qquad \text{for all } y \in \mathbb{R}^n, \alpha \in \mathbb{R}.$$

We also have two crucial properties which follow directly from (2.12):

$$(2.14) \qquad G(py) = pG(y), \qquad S(py) = pS(y) \quad \text{for all } p \in \Sigma_n, y \in \mathcal{R}^\circ.$$

Equations (2.3) and (2.5) may be rewritten in this setting as

$$(2.3)' \qquad dY(t) = (G(Y(t)) + \gamma v) \, dt + S(Y(t)) \, dW(t), \qquad Y(0) = \bar{y},$$

$(2.5)'$ $\quad d\left(\sum_{i=1}^{n} Y_i(t)\right) = n\gamma\,dt + v^t S(Y(t))\,dW(t),$

respectively. Now define the process

$$\widetilde{Y}(t) := Y(t) - \left(\frac{1}{n}\sum_{i=1}^{n} Y_i(t)\right)v, \qquad 0 \le t < \infty,$$

which lives in the subspace $\Pi := \{(y_1,\ldots,y_n) \in \mathbb{R}^n \,|\, y_1 + \cdots + y_n = 0\}$ of $\mathbb{R}^n$ with normal vector $v$. From $(2.3)'$ and $(2.5)'$, we have

$$d\widetilde{Y}(t) = G(Y(t))\,dt + \widetilde{S}(Y(t))\,dW(t), \qquad \widetilde{Y}(0) = \tilde{y},$$

where $\tilde{y} := \bar{y} - (\frac{1}{n}\sum_{i=1}^{n} \bar{y}_i)v$ and $\widetilde{S}(y) := S(y) - \frac{1}{n}vv^t S(y)$ for all $y \in \mathbb{R}^n$. In fact,

$(2.15)$ $\qquad d\widetilde{Y}(t) = G(\widetilde{Y}(t))\,dt + \widetilde{S}(\widetilde{Y}(t))\,dW(t), \qquad \widetilde{Y}(0) = \tilde{y},$

because of $(2.13)$. We note that if $x \in \Pi$, then for any $y \in \mathbb{R}^n$,

$$x^t \widetilde{S}(y)x = x^t S(y)x - \frac{1}{n}x^t vv^t S(y)x = x^t S(y)x \ge \sigma_{\min}\|x\|^2,$$

where we have set $\sigma_{\min} := \min\{\sigma_1,\ldots,\sigma_n\} > 0$ and used the fact that $x^t v \equiv x \cdot v = 0$. This means that the covariance matrix in $(2.15)$ is uniformly nondegenerate when restricted to the subspace $\Pi$. In particular, the theory of Bass and Pardoux [1] once again shows that the $\Pi$-valued solution $\widetilde{Y}(\cdot)$ of $(2.15)$ is unique in the sense of the probability law.

We now claim that

$(2.16)$ $\qquad\qquad y \cdot G(y) \le c\|y\| \qquad$ holds for all $y \in \Pi$,

where $c < 0$ is a constant depending only on $n$ and $g_1,\ldots,g_n$. Indeed, fix $y \in \Pi$. There exists $p \in \Sigma_n$ such that $y \in \mathcal{R}_p$, so

$(2.17)$
$$\begin{aligned}
y \cdot G(y) &= \sum_{i=1}^{n} y_i g_{p^{-1}(i)} = \sum_{k=1}^{n} y_{p(k)} g_k \\
&= y_{p(n)} \sum_{m=1}^{n} g_m + \sum_{k=1}^{n-1}\left(y_{p(k)} - y_{p(k+1)}\right)\sum_{m=1}^{k} g_m,
\end{aligned}$$

where the final equality follows by summation by parts. From $(1.5)$, we have $\sum_{m=1}^{n} g_m = 0$ and $\sum_{m=1}^{k} g_m < 0$ for all $k = 1,\ldots,n-1$, and $(2.11)$ gives $y_{p(k)} - y_{p(k+1)} \ge 0$ for the same range of $k$. Set $c = n^{-1/2}\max_{1 \le k \le n-1}\{\sum_{m=1}^{k} g_m\} < 0$, and note that $(2.11)$ and the fact that $y \in \Pi$ imply that $y_{p(1)} \ge 0$, $y_{p(n)} \le 0$ and $y_{p(1)} \geq y_i$ for all $i = 1, \ldots, n$. In particular, $\|y\|^2 \leq n \max\{y_{p(1)}^2, y_{p(n)}^2\} \leq n(y_{p(1)} - y_{p(n)})^2$. Finally, (2.17) gives

$$y \cdot G(y) = \sum_{k=1}^{n-1} (y_{p(k)} - y_{p(k+1)}) \sum_{m=1}^{k} g_m \leq c\sqrt{n} \sum_{k=1}^{n-1} (y_{p(k)} - y_{p(k+1)})$$
$$= c\sqrt{n}(y_{p(1)} - y_{p(n)}) \leq c\|y\|.$$

In the Appendix, it is shown that (2.16) implies that the process

$$
\begin{align}
(2.18) \qquad & \widetilde{Y}(\cdot) \text{ is recurrent with respect to } B \cap \Pi, \\
& \text{for some ball } B \subset \mathbb{R}^n \text{ centered at } 0.
\end{align}
$$

Theorem 5.1 on page 121 of [11] guarantees that the process $\widetilde{Y}(\cdot)$ of (2.15) admits a stationary distribution $\mu$, such that for any bounded, measurable function $f : \Pi \to \mathbb{R}$ we have

$$(2.19) \qquad \lim_{T \to \infty} \frac{1}{T} \int_0^T f(\widetilde{Y}(t)) \, dt = \int_\Pi f(y) \, d\mu(y) \qquad \text{a.s.}$$

Moreover, $\mu$ is a probability measure on $\Pi$ which does not depend on the initial value $\widetilde{Y}(0) = \tilde{y}$. Setting $\bar{y} = 0 = \tilde{y}$, fix some $p \in \Sigma_n$ and apply it to (2.15) to obtain

$$d(p\widetilde{Y}(t)) = pG(\widetilde{Y}(t)) \, dt + p\widetilde{S}(\widetilde{Y}(t)) \, dW(t), \qquad \widetilde{Y}(0) = 0,$$

which, in view of (2.14), may be rewritten as

$$(2.20) \quad d(p\widetilde{Y}(t)) = G_1(p\widetilde{Y}(t)) \, dt + \widetilde{S}_1(p\widetilde{Y}(t)) \, dW(t), \qquad p\widetilde{Y}(0) = 0,$$

where $G_1(y) = G(y)$ and $\widetilde{S}_1(y) = \widetilde{S}(y)$ for all $y$ not in the exceptional set $\mathcal{R}^e$. It is argued in the Appendix that (2.15), (2.20) imply that the processes

$$(2.21) \qquad \widetilde{Y}(\cdot) \text{ and } p\widetilde{Y}(\cdot) \text{ have the same stationary distribution } \mu.$$

Since $p$ is arbitrary, it follows that $\mu$ is invariant under the action of $\Sigma_n$. From (2.19) with $f(y) := \mathbb{1}_{\mathcal{R}_p}(y)$, we obtain

$$\lim_{T \to \infty} \frac{1}{T} \int_0^T \mathbb{1}_{\mathcal{R}_p}(\widetilde{Y}(t)) \, dt = \int_\Pi \mathbb{1}_{\mathcal{R}_p}(y) \, d\mu(y) = \int_\Pi \mathbb{1}_{\mathcal{R}_1}(y) \, d\mu(y) \qquad \text{a.s.}$$

By the remark preceding (2.13), we may replace $\widetilde{Y}(t)$ by $Y(t)$ in the above equation to conclude that the a.s.-limiting value of $T^{-1} \int_0^T \mathbb{1}_{\mathcal{R}_p}(Y(t)) \, dt$ is independent of $p$. Summing over all $p \in \Sigma_n$, we find that

$$(2.22) \qquad \lim_{T \to \infty} \frac{1}{T} \int_0^T \mathbb{1}_{\mathcal{R}_p}(Y(t)) \, dt = \frac{1}{n!} \qquad \text{for all } p \in \Sigma_n, \text{ a.s.}$$

For fixed $i$ and $k$, (2.10) now follows by summing (2.22) over the $(n-1)!$ permutations $p \in \Sigma_n$ satisfying $p(k) = i$. $\quad \square$

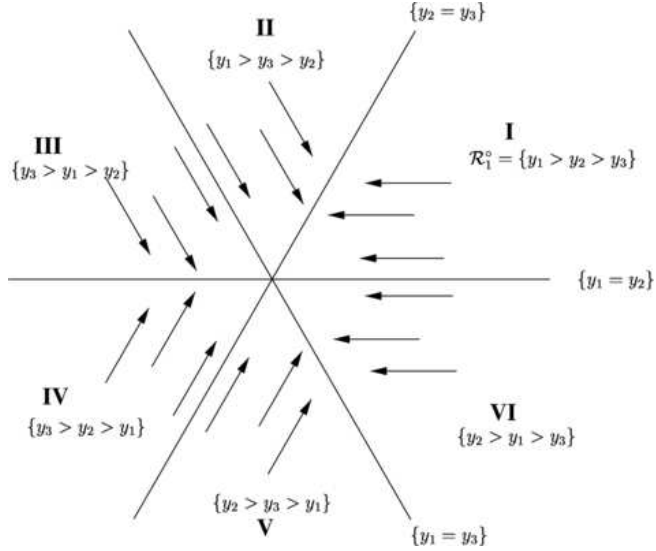

FIG. 1.   *Projection of the Atlas model with $n = 3$ onto the subspace* $\Pi$.

The above proof shows that any given ranking of the stocks in a first-order model occurs roughly $1/n!$ of the time. This does not imply that the rank changes occur with roughly the same frequency at all scales. For example, consider an Atlas model with $n = 3$ and constant volatilities across ranks. The plane $\Pi = \{(y_1, y_2, y_3) | y_1 + y_2 + y_3 = 0\}$ is represented in Figure 1. The vectors shown are values for $G(y)$ for various $y \in \Pi$. The function $G(\cdot)$ is constant within each of the six wedges; furthermore, changes in rank occur when the process $\widetilde{Y}(\cdot)$ hits the exceptional set $\mathcal{R}^e$, which is the union of the three lines shown. It is clear from the direction of the vectors that changes of rank will be likely to occur much more frequently between the bottom two stocks than the top two stocks. That is, changes in $\widetilde{Y}(\cdot)$ between the three pairs of regions labelled (I, II), (III, IV) and (V, VI) occur more frequently than between (II, III), (IV, V) and (VI, I). Of course, this does not hold in the general first-order model. The above proof also reveals that the rank-ordered process $\sum_{p \in \Sigma_n} \mathbb{1}_{\mathcal{R}_p}(Y(\cdot)) p^{-1} Y(\cdot)$ is a reflected Brownian motion in the polyhedral region $\mathcal{R}_1$ with constant drift equal to $(g_1 + \gamma, \ldots, g_n + \gamma)$, covariance matrix given by $\mathrm{diag}(\sigma_1, \ldots, \sigma_n)$, and normal reflection on the boundary.

**3. Ranked capitalization processes.**   Having constructed the solution $Y(\cdot) = (Y_1(\cdot), \ldots, Y_n(\cdot))$ of the stochastic differential system (2.3), let us now look at the processes

$$(3.1) \qquad Z_k(t) := \sum_{i=1}^{n} \mathbb{1}_{\mathcal{Q}_k^{(i)}}(Y(t)) \cdot Y_i(t), \qquad 0 \le t < \infty,$$

for $k = 1, \ldots, n$. These are the log-capitalizations of the various companies listed according to their rank, so that

$$(3.2) \qquad X_{(k)}(t) = e^{Z_k(t)}, \qquad \mu_{(k)}(t) := \frac{X_{(k)}(t)}{X_1(t) + \cdots + X_n(t)}$$

represent, respectively, the absolute and relative capitalizations of the company ranked $k$th at time $t$, in accordance with (1.3) and (2.8). Denoting by $\Lambda^{k,k+1}(\cdot) := \Lambda_{Z_k - Z_{k+1}}(\cdot)$ the local time accumulated *at the origin* by the nonnegative semimartingale $Z_k(\cdot) - Z_{k+1}(\cdot)$ up to calendar time $t$, and setting

$$\Lambda^{0,1}(\cdot) \equiv 0, \qquad \Lambda^{n,n+1}(\cdot) \equiv 0,$$

we obtain the dynamics for the processes in (3.1) in the form

$$dZ_k(t) = \sum_{i=1}^n \mathbb{1}_{Q_k^{(i)}(Y(t))} \cdot dY_i(t) + \tfrac{1}{2}[d\Lambda^{k,k+1}(t) - d\Lambda^{k-1,k}(t)],$$

or equivalently

$$(3.3) \qquad Z_k(t) = Z_k(0) + (g_k + \gamma)t + \sigma_k B_k(t) + \tfrac{1}{2}[\Lambda^{k,k+1}(t) - \Lambda^{k-1,k}(t)],$$
$$0 \le t < \infty.$$

We have used the equations of (2.3) and the notation of (2.5), and have applied the generalized Itô rule for convex functions of semimartingales from Section 3.7 in [8], in a manner similar to the derivations in Chapter 4 of [2]. (These derivations require that the processes $Y_1, \ldots, Y_n$ be *pathwise mutually nondegenerate*, as in Definition 4.1.2 of [2]; however, this follows from an application of the Girsanov theorem, which is justified by the uniform nondegeneracy of the variance structure and boundedness of the drift coefficients.)

In conjunction now with (2.7), the dynamics (3.3) yield the strong law of large numbers

$$(3.4) \qquad \lim_{T \to \infty} \frac{1}{T}[\Lambda^{k-1,k}(T) - \Lambda^{k,k+1}(T)] = 2g_k \qquad \text{a.s.}$$

for every $k = 1, \ldots, n$. Taking $k = 1$, this means that the limit

$$\lambda_{1,2} := \lim_{T \to \infty} \frac{1}{T}\Lambda^{1,2}(T)$$

exists a.s., and that $\lambda_{1,2} = -2g_1$. Arguing by induction, we see that all limits

$$(3.5) \qquad \lambda_{k,k+1} := \lim_{T \to \infty} \frac{1}{T}\Lambda^{k,k+1}(T)$$

exist a.s. and satisfy

$$(3.6) \qquad \lambda_{k-1,k} - \lambda_{k,k+1} = 2g_k \qquad \text{for } k = 1, \ldots, n$$

(of course, $\lambda_{0,1} = \lambda_{n,n+1} = 0$). In other words, the quantities of (3.5) are given as

$$(3.7) \qquad \lambda_{k,k+1} = -2(g_1 + \cdots + g_k) > 0 \qquad \text{for } k = 1, \ldots, n-1.$$

Observe now from (3.3) the decomposition

$$(3.8) \qquad Z_k(t) - Z_{k+1}(t) = Z_k(0) - Z_{k+1}(0) + \Theta_k(t) + \Lambda_{Z_k - Z_{k+1}}(t)$$

for the nonnegative semimartingale $Z_k(\cdot) - Z_{k+1}(\cdot)$. We are using here the notation

$$(3.9) \quad \Theta_k(t) := (g_k - g_{k+1})t - \frac{1}{2}[\Lambda^{k-1,k}(t) + \Lambda^{k+1,k+2}(t)] + s_k \cdot \widetilde{W}^{(k)}(t),$$

where

$$(3.10) \quad s_k := \sqrt{\sigma_k^2 + \sigma_{k+1}^2} \quad \text{and} \quad \widetilde{W}^{(k)}(t) := \frac{1}{s_k}(\sigma_k B_k(t) - \sigma_{k+1} B_{k+1}(t))$$

is standard Brownian motion. This decomposition (3.8) shows that $Z_k(\cdot) - Z_{k+1}(\cdot)$ is the *reflection at the origin* of the semimartingale $\Theta_k(\cdot)$ in (3.9). Now the bounded variation part of the semimartingale $\Theta_k(\cdot)$ is of the form

$$\left\{ (g_k - g_{k+1}) - \frac{1}{2t}[\Lambda^{k-1,k}(t) + \Lambda^{k+1,k+2}(t)] \right\} t = -(\lambda_{k,k+1} + o(1))t$$

as $t \to \infty$, thanks to (3.5) and (3.6). Thus $\Theta_k(\cdot)$ behaves asymptotically as Brownian motion with negative drift $-\lambda_{k,k+1}$.

## 4. Stability of capital distribution.

Let us look now at the ergodic behavior, as $t \to \infty$, of the nonnegative process

$$(4.1) \qquad \Xi_k(t) := \log\left(\frac{\mu_{(k)}(t)}{\mu_{(k+1)}(t)}\right) = Z_k(t) - Z_{k+1}(t) = \Xi_k(0) + \Theta_k(t) + \Lambda_{\Xi_k}(t),$$
$$0 \le t < \infty,$$

in (3.8): namely, the reflection at the origin of the semimartingale $\Theta_k(\cdot)$ of (3.9). As we have remarked, this process $\Xi_k(\cdot)$ behaves asymptotically as Brownian motion with negative drift $-\lambda_{k,k+1}$, reflected at the origin. Therefore,

$$(4.2) \qquad \lim_{t \to \infty} \log\left(\frac{\mu_{(k)}(t)}{\mu_{(k+1)}(t)}\right) = \lim_{t \to \infty} \Xi_k(t) = \xi_k \qquad \text{in distribution.}$$

Here, for each $k = 1, \ldots, n-1$ the random variable $\xi_k$ has an *exponential distribution with parameter*

$$(4.3) \qquad r_k := \frac{2\lambda_{k,k+1}}{s_k^2} = -\frac{4(g_1 + \cdots + g_k)}{\sigma_k^2 + \sigma_{k+1}^2} > 0$$
$$\text{that is, } \mathbb{P}(\xi_k > x) = e^{-r_k x} \ \forall \, x \ge 0.$$

This leads to the *asymptotic Pareto distribution*

$$(4.4) \qquad \lim_{t \to \infty} \mathbb{P}\left[\frac{\mu_{(k)}(t)}{\mu_{(k+1)}(t)} > y\right] = y^{-r_k} = \mathbb{P}(\xi_k > \log y) \qquad \forall\, y \geq 1$$

for the ratios of successively ranked capitalizations, which is frequently observed in practice; see Chapter 5 of [2], in particular Figure 5.1 on page 95 and the discussion on page 102.

We also obtain for every $k = 1, \ldots, n-1$ the strong law of large numbers

$$(4.5) \qquad \lim_{T \to \infty} \frac{1}{T} \int_0^T g(\Xi_k(t))\, dt = \lim_{T \to \infty} \frac{1}{T} \int_0^T g\left(\log \frac{\mu_{(k)}(t)}{\mu_{(k+1)}(t)}\right) dt$$
$$= \mathbb{E}[g(\xi_k)] \qquad \text{a.s.}$$

for every measurable function $g : [0, \infty) \to \mathbb{R}$ with $\int_0^\infty |g(x)|e^{-r_k x}\, dx < \infty$; see [10], Theorem 3.1, [5], Section 23 or [11], Theorem 5.1 on page 121.

In fact, we can ascertain a little more generally that we have *in distribution*:

$$(4.6) \qquad \lim_{t \to \infty} (\Xi_1(t), \ldots, \Xi_{n-1}(t)) = \lim_{t \to \infty} \left(\log \frac{\mu_{(1)}(t)}{\mu_{(2)}(t)}, \ldots, \log \frac{\mu_{(n-1)}(t)}{\mu_{(n)}(t)}\right)$$
$$= (\xi_1, \ldots, \xi_{n-1}).$$

This follows from the very detailed analysis of what Harrison and Williams [6, 7] and Williams [12] call "*reflected* (or regulated) *Brownian motions*" (RBMs, for short) in polyhedral domains; see, in particular, Sections 4–8 of [7] which are of particular relevance to our setting here. Then Theorem 3.1 of [10] guarantees again that a strong law of large numbers

$$(4.7) \qquad \lim_{T \to \infty} \frac{1}{T} \int_0^T g(\Xi_1(t), \ldots, \Xi_{n-1}(t))\, dt = \mathbb{E}[g(\xi_1, \ldots, \xi_{n-1})]$$

holds a.s., for every bounded, measurable $g : [0, \infty)^{n-1} \to \mathbb{R}$.

REMARK. As discussed in [6, 7] and [13], the joint distribution of the random vector $(\xi_1, \ldots, \xi_{n-1})$ has a probability density function $\mathfrak{p} : [0, \infty)^{n-1} \to [0, \infty)$ with $\int_{[0,\infty)^{n-1}} \mathfrak{p}(x)\, dx = 1$ that satisfies a certain integral equation (the *basic adjoint relation* of (3.2) in [13] or (BAR) on page 103 in [7]). This equation involves the second-order diffusion operator

$$L\mathfrak{p}(x) := \frac{1}{2} \sum_{k=1}^{n-1} \sum_{\ell=1}^{n-1} A_{k\ell} \frac{\partial^2 \mathfrak{p}(x)}{\partial x_k\, \partial x_\ell} + \sum_{k=1}^{n-1} \lambda_{k,k+1} \frac{\partial \mathfrak{p}(x)}{\partial x_k},$$

and specifies appropriate boundary conditions on the faces of the orthant $[0, \infty)^{n-1}$. Here $A = \{A_{k\ell}\} = \Sigma\Sigma'$, where we have set

$$\Sigma = \begin{pmatrix} \sigma_1 & -\sigma_2 & 0 & \cdots & 0 & 0 \\ 0 & \sigma_2 & -\sigma_3 & \ldots & 0 & 0 \\ \vdots & \vdots & \vdots & \ddots & \vdots & \vdots \\ 0 & 0 & 0 & \cdots & \sigma_{n-1} & -\sigma_n \end{pmatrix}$$

for the $(n-1) \times n$ volatility matrix for the *multidimensional Brownian motion*

$$\widehat{\Xi}_k(t) = \widehat{\Xi}_k(0) - \lambda_{k,k+1}t + \sigma_k B_k(t) - \sigma_{k+1} B_{k+1}(t) + \Lambda_{\widehat{\Xi}_k}(t), \qquad 0 \le t < \infty,$$

for $k = 1, \dots, n-1$, *with normal reflection on each of the faces of the orthant.* In particular,

$$A = \begin{pmatrix} \sigma_1^2 + \sigma_2^2 & -\sigma_2^2 & 0 & \cdots & 0 & 0 \\ -\sigma_2^2 & \sigma_2^2 + \sigma_3^2 & -\sigma_3^2 & \cdots & 0 & 0 \\ 0 & -\sigma_3^2 & \sigma_3^2 + \sigma_4^2 & \cdots & 0 & 0 \\ \vdots & \vdots & \vdots & \ddots & \vdots & \vdots \\ 0 & 0 & 0 & \cdots & -\sigma_{n-1}^2 & \sigma_{n-1}^2 + \sigma_n^2 \end{pmatrix}.$$

The structure of this process has slightly more complicated volatilities than the setting of [7]; satisfies their "stability condition" $\lambda_{1,2} > 0, \dots, \lambda_{n-1,n} > 0$ which ensures positive recurrence; but fails to satisfy their "skew-symmetry condition," that makes $\mathfrak{p}(\cdot)$ the product of the exponential densities in (4.3). It is highly unlikely that the "basic adjoint relation" which characterizes $\mathfrak{p}(\cdot)$ can be solved in closed form; as a result, we know only the one-dimensional marginals of the density $\mathfrak{p}(\cdot)$, not the density itself.

By the Skorohod representation, one can construct now, possibly on an enlarged probability space, copies $(\widetilde{\xi}_1, \dots, \widetilde{\xi}_{n-1})$ of the random variables $(\xi_1, \dots, \xi_{n-1})$, as well as copies $(\widetilde{\mu}_1(\cdot), \dots, \widetilde{\mu}_{n-1}(\cdot))$ of the processes $(\mu_1(\cdot), \dots, \mu_{n-1}(\cdot))$, such that

$$\lim_{t \to \infty} \left( \log \frac{\widetilde{\mu}_{(1)}(t)}{\widetilde{\mu}_{(2)}(t)}, \dots, \log \frac{\widetilde{\mu}_{(n-1)}(t)}{\widetilde{\mu}_{(n)}(t)} \right) = (\widetilde{\xi}_1, \dots, \widetilde{\xi}_{n-1})$$

holds almost surely instead of just in distribution as in (4.6). In particular,

$$\lim_{t \to \infty} \left( \frac{\widetilde{\mu}_{(k)}(t)}{\widetilde{\mu}_{(n)}(t)} \right) = e^{\widetilde{\xi}_{n-1} + \cdots + \widetilde{\xi}_k}, \qquad k = 1, \dots, n-1,$$

holds almost surely, and gives

$$\lim_{t \to \infty} \left( \frac{1 - \widetilde{\mu}_{(n)}(t)}{\widetilde{\mu}_{(n)}(t)} \right) = e^{\widetilde{\xi}_{n-1}} + e^{\widetilde{\xi}_{n-1} + \widetilde{\xi}_{n-2}} + \cdots + e^{\widetilde{\xi}_{n-1} + \cdots + \widetilde{\xi}_1} \qquad \text{a.s.}$$

which then leads to

$$(4.8) \qquad \lim_{t \to \infty} \widetilde{\mu}_{(k)}(t) = \widetilde{M}_k := \frac{e^{\widetilde{\xi}_{n-1} + \cdots + \widetilde{\xi}_k}}{1 + e^{\widetilde{\xi}_{n-1}} + \cdots + e^{\widetilde{\xi}_{n-1} + \cdots + \widetilde{\xi}_1}} \qquad \text{a.s.}$$

for every $k = 1, \dots, n$. The understanding here is that the "empty summation," which occurs in the numerator when $k = n$, is taken to be equal to zero.

Consider now the vector of random variables $(M_1, \ldots, M_n)$ defined by

$$(4.9) \quad M_n := (1 + e^{\xi_{n-1}} + \cdots + e^{\xi_{n-1} + \cdots + \xi_1})^{-1} \quad \text{and} \quad M_k := M_n \cdot e^{\xi_{n-1} + \cdots + \xi_k}$$

for $k = 1, \ldots, n-1$, by analogy with (4.8). From the equation (4.8) and the discussion that precedes it, we can recast (4.6) as

$$(4.10) \quad \lim_{t \to \infty} (\mu_{(1)}(t), \ldots, \mu_{(n)}(t)) = (M_1, \ldots, M_n) \qquad \text{in distribution.}$$

*The quantities of (4.9) are the long-term relative weights of the various stocks in this market*, represented by their ranks—from the largest (namely, $M_1$) down to the smallest (namely, $M_n$). We also have from (4.7) the strong law of large numbers

$$(4.11) \quad \lim_{T \to \infty} \frac{1}{T} \int_0^T f(\mu_{(1)}(t), \ldots, \mu_{(n)}(t)) \, dt = \mathbb{E}[f(M_1, \ldots, M_n)] \qquad \text{a.s.}$$

for every bounded and measurable $f : \Delta^n \to \mathbb{R}$, where $\Delta^n = \{(x_1, \ldots, x_n) \in \mathbb{R}^n | x_1 \geq 0, \ldots, x_n \geq 0, \sum_{j=1}^n x_j = 1\}$.

4.1. *The certainty-equivalent approximation.* The random vector $M = (M_1, \ldots, M_n)$ of (4.9), (4.10) is hard to come to grips with: as we have already remarked, we do not know much about the joint distribution of the random variables $(\xi_1, \ldots, \xi_{n-1})$ which determine it. In fact, we only know the one-dimensional marginal distributions of the random variables $\xi_1, \ldots, \xi_{n-1}$ individually—namely, the exponentials of (4.3). We look then at a particularly convenient approximation of the vector $M = (M_1, \ldots, M_n)$, obtained by replacing the random variables $\xi_k$ in (4.9) by their expected values

$$(4.12) \quad \rho_k := \mathbb{E}(\xi_k) = 1/r_k = -\frac{\sigma_k^2 + \sigma_{k+1}^2}{4(g_1 + \cdots + g_k)} > 0, \qquad k = 1, \ldots, n-1,$$

namely,

$$(4.13) \quad \begin{aligned} M_n^{CE} &:= (1 + e^{\rho_{n-1}} + \cdots + e^{\rho_{n-1} + \cdots + \rho_1})^{-1} \quad \text{and} \\ M_k^{CE} &:= M_n^{CE} \cdot e^{\rho_{n-1} + \cdots + \rho_k} \end{aligned}$$

for $k = 1, \ldots, n-1$. We call the resulting probability vector $M^{CE} = (M_1^{CE}, \ldots, M_n^{CE})$ the *certainty-equivalent approximation* of the long-term relative capitalizations $(M_1, \ldots, M_n)$ in (4.9), (4.10).

EXAMPLE 4.1. For the Atlas model of (1.7) with constant variances

$$(4.14) \qquad \sigma_k^2 = \sigma^2 > 0 \qquad \text{for } k = 1, \ldots, n,$$

the quantities of (4.12) take the form $\rho_k = \alpha/k, k = 1, \ldots, n-1$ with $\alpha := \sigma^2/2g$, and we have a further approximation for the certainty-equivalent quantities of (4.13), namely,

$$(4.15) \qquad M_k^{CE} \sim \frac{k^{-\alpha}}{\sum_{j=1}^n j^{-\alpha}}, \qquad k = 1, \ldots, n.$$

EXAMPLE 4.2. For the generalized Atlas model of (1.7) with linearly growing variances

$$(4.16) \qquad \sigma_k^2 = \sigma^2 + ks^2, \qquad k = 1, \ldots, n \text{ for some } \sigma^2 > 0, s^2 \geq 0,$$

we get $\rho_k = 2\beta + \frac{\alpha+\beta}{k}, k = 1, \ldots, n-1$ with $\beta := s^2/4g$, and the certainty-equivalents of (4.13) are now approximated as

$$(4.17) \qquad M_k^{CE} \sim \frac{k^{-(\alpha+\beta)} e^{-2\beta k}}{\sum_{j=1}^n j^{-(\alpha+\beta)} e^{-2\beta j}}, \qquad k = 1, \ldots, n.$$

**5. Portfolios and their growth rates.** Let us consider now investing in the market of (1.1), (1.6)—equivalently modeled by (2.3) and (2.4)—according to a portfolio rule $\pi = (\pi_1, \ldots, \pi_n)$. This is a process adapted to the natural filtration $\mathbb{F} = \{\mathcal{F}(t)\}_{0 \leq t < \infty}$ of the stock-prices $\mathcal{F}(t) := \sigma(X(s), \ 0 \leq s \leq t)$, which satisfies $\pi_1(t) \geq 0, \ldots, \pi_n(t) \geq 0$ and $\sum_{i=1}^n \pi_i(t) = 1$ for all $0 \leq t < \infty$. The interpretation is that $\pi_i(t)$ represents the proportion of the portfolio's wealth $Z^\pi(t)$ that is invested at time $t$ in the $i$th stock, so that

$$(5.1) \qquad \begin{aligned} \frac{dZ^\pi(t)}{Z^\pi(t)} &= \sum_{i=1}^n \pi_i(t) \cdot \frac{dX_i(t)}{X_i(t)} \\ &= \sum_{i=1}^n \pi_i(t) \cdot \left[ \left( \gamma_i(t) + \frac{\sigma_i^2(t)}{2} \right) dt + \sigma_i(t) \, dW_i(t) \right], \end{aligned}$$

and $Z^\pi(0) = z > 0$ is the initial capital. The quantities $b_i(t) := \gamma_i(t) + (\sigma_i^2(t)/2)$ for $i = 1, \ldots, n$ appearing in (5.1) are the *rates of return* of the individual stocks.

As shown in [2], Chapter 1, an application of Itô's rule casts (5.1) in the equivalent form

$$(5.2) \qquad \log Z^\pi(T) = \log z + \int_0^T \gamma^\pi(t) \, dt + \sum_{i=1}^n \int_0^T \pi_i(t) \sigma_i(t) \, dW_i(t),$$
$$0 \leq t < \infty.$$

Here the quantities

$$(5.3) \quad \gamma^\pi(t) := \sum_{i=1}^n \pi_i(t) \gamma_i(t) + \gamma_*^\pi(t), \qquad \gamma_*^\pi(t) := \frac{1}{2} \sum_{i=1}^n \pi_i(t)(1 - \pi_i(t)) \sigma_i^2(t)$$

denote, respectively, the *growth rate* and the *excess growth rate* of the portfolio. The appellation is justified by the a.s. equality

$$(5.4) \qquad \lim_{T \to \infty} \frac{1}{T} \Big( \log Z^\pi(T) - \int_0^T \gamma^\pi(t) \, dt \Big) = 0;$$

this is a consequence of (5.2), the boundedness of $\pi_i(\cdot)$, $\sigma_i(\cdot)$ and the strong law of large numbers for Brownian motion. The *rate of return* and the *variance* of the portfolio $\pi(\cdot)$ are, respectively, the quantities

$$(5.5) \quad b^\pi(t) := \sum_{i=1}^n \pi_i(t)(\gamma_i(t) + (\sigma_i^2(t)/2)) \quad \text{and} \quad (\sigma^\pi(t))^2 := \sum_{i=1}^n (\pi_i(t))^2 \sigma_i^2(t).$$

We shall denote by

$$(5.6) \quad G^\pi(n) := \lim_{T \to \infty} \frac{1}{T} \int_0^T \gamma^\pi(t) \, dt, \qquad G_*^\pi(n) := \lim_{T \to \infty} \frac{1}{T} \int_0^T \gamma_*^\pi(t) \, dt$$

the long-term averages of the growth rate and of the excess growth rate, respectively, whenever these limits exist a.s. In fact, when the first limit in (5.6) exists, it is clear from (5.4) that

$$G^\pi(n) = \lim_{T \to \infty} \frac{1}{T} \log Z^\pi(T)$$

will also hold a.s. We have parametrized the quantities of (5.6) by the market size $n$ because we shall also be interested in the large-market behavior, as the number of equities tends to infinity, of some simple and consistently defined families $\Pi = \{\pi^{(n)}\}_{n \in \mathbb{N}}$ of portfolio rules:

$$(5.7) \qquad \Gamma^\Pi := \lim_{n \to \infty} G^{\pi^{(n)}}(n), \qquad \Gamma_*^\Pi := \lim_{n \to \infty} G_*^{\pi^{(n)}}(n).$$

We shall study, and then compare to each other, the quantities of (5.6), (5.7) for a few simple but important and relatively easy-to-implement investment rules, for which the limits indicated in (5.7) exist.

EXAMPLE 5.1. The market portfolio rule $\mu(\cdot) = (\mu_1(\cdot), \dots, \mu_n(\cdot))$ has already been introduced in (2.8), namely,

$$(5.8) \qquad \mu_i(t) := \frac{X_i(t)}{X_1(t) + \dots + X_n(t)}, \qquad i = 1, \dots, n.$$

It invests in each company in proportion to its relative market capitalization, and yields a wealth $Z^\mu(\cdot) = zX(\cdot)$ that reflects the entire market capitalization, in proportion of course to the initial investment $z > 0$. For this market portfolio, and with the notation of (3.2) for the reverse order-statistics, the quantities of (5.3) become

$$(5.9) \qquad \begin{aligned} \gamma^\mu(t) &= \gamma + \sum_{i=1}^n g_k \mu_{(k)}(t) + \gamma_*^\mu(t), \\ \gamma_*^\mu(t) &= \tfrac{1}{2} \sum_{k=1}^n \sigma_k^2 \cdot \mu_{(k)}(t)(1 - \mu_{(k)}(t)) \end{aligned}$$

in the context of the model of (1.1), (1.6). Also, we know from Remark 2.1 that

$$(5.10) \qquad G^\mu(n) := \lim_{T \to \infty} \frac{1}{T} \int_0^T \gamma^\mu(t)\,dt = \lim_{T \to \infty} \frac{1}{T} \log Z^\mu(T) = \gamma \qquad \text{a.s.}$$

so that (5.9), (4.11) then imply

$$(5.11) \qquad \begin{aligned} G_*^\mu(n) &= \lim_{T \to \infty} \frac{1}{T} \int_0^T \gamma_*^\mu(t)\,dt \\ &= \frac{1}{2} \sum_{k=1}^n \sigma_k^2 \cdot \mathbb{E}(M_k(1 - M_k)) = -\sum_{k=1}^n g_k \cdot \mathbb{E}(M_k). \end{aligned}$$

Here $(M_1, \ldots, M_n)$ is the vector of long-term ranked market weights of (4.9), (4.10).

For the (generalized) Atlas model of (1.7), these formulae give

$$(5.12) \qquad G^\mu(n) = g, \qquad G_*^\mu(n) = g(1 - n \cdot \mathbb{E}(M_n)).$$

EXAMPLE 5.2.   The equally-weighted portfolio rule $\eta(\cdot) = (\eta_1(\cdot), \ldots, \eta_n(\cdot))$ assigns equal weights

$$(5.13) \qquad \eta_i(t) := \frac{1}{n}, \qquad i = 1, \ldots, n,$$

to all stocks at all times. Clearly,

$$(5.14) \qquad \gamma^\eta(t) = \gamma + \gamma_*^\eta(t), \qquad \gamma_*^\eta(t) = \frac{n-1}{2n^2} \sum_{k=1}^n \sigma_k^2.$$

For the variance structure of (4.16) and with $\gamma = g > 0$, these imply

$$(5.15) \qquad G_*^\eta(n) = \frac{n-1}{2n}\left(\sigma^2 + s^2 \frac{n+1}{2}\right), \qquad G^\eta(n) = g + G_*^\eta(n).$$

In particular, for the constant-variance case of (4.14) we get

$$(5.16) \qquad \Gamma_*^\eta = \frac{\sigma^2}{2}, \qquad \Gamma^\eta = g + \frac{\sigma^2}{2} \qquad \text{when } s = 0,$$

whereas $\Gamma_*^\eta = \Gamma^\eta = \infty$ when $s > 0$ in (4.16).

EXAMPLE 5.3.  The diversity-weighted portfolio rule $\vartheta^{(p)}(\cdot) = (\vartheta_1^{(p)}(\cdot), \ldots, \vartheta_n^{(p)}(\cdot))$ is given in terms of the market portfolio as

$$(5.17) \qquad \vartheta_i^{(p)}(t) := \frac{(\mu_i(t))^p}{(\mu_1(t))^p + \cdots + (\mu_n(t))^p}, \qquad i = 1, \ldots, n,$$

for some fixed number $p \in (0,1)$. This portfolio has been studied already by Fernholz [2] and Fernholz, Karatzas and Kardaras [4]; in particular, we know from these sources that

$$(5.18) \qquad \log\left(\frac{Z^{\vartheta^{(p)}}(T)}{Z^\mu(T)}\right) = \log\left(\frac{D(\mu(T))}{D(\mu(0))}\right) + (1-p) \int_0^T \gamma_*^{\vartheta^{(p)}}(t)\,dt,$$
$$0 \le T < \infty,$$

holds a.s. for every $T \in (0, \infty)$, with the notation $D(x) := (\sum_{i=1}^n x_i^p)^{1/p}$. But (5.3) gives

$$\gamma^{\vartheta^{(p)}}(t) - \gamma_*^{\vartheta^{(p)}}(t) = \gamma + \frac{\sum_{k=1}^n g_k(\mu_{(k)}(t))^p}{\sum_{k=1}^n (\mu_{(k)}(t))^p},$$

with the help of which (5.18) reads

$$\frac{1}{T}\left(\log Z^{\vartheta^{(p)}}(T) - \int_0^T \gamma^{\vartheta^{(p)}}(t)\,dt\right)$$
$$= \frac{1}{T}(\log Z^\mu(T) - \gamma) - \frac{1}{T}\int_0^T \frac{\sum_{k=1}^n g_k(\mu_{(k)}(t))^p}{\sum_{k=1}^n (\mu_{(k)}(t))^p}\,dt$$
$$+ \frac{1}{T}\log\left(\frac{D(\mu(T))}{D(\mu(0))}\right) - \frac{p}{T}\int_0^T \gamma_*^{\vartheta^{(p)}}(t)\,dt.$$

Sending $T$ to infinity in this expression and recalling (5.4), (5.10) and (4.11), we deduce that the a.s. limits of (5.6) exist for the diversity-weighted portfolio, and are given by

$$(5.19) \qquad G_*^{\vartheta^{(p)}}(n) = -\frac{1}{p} \cdot \mathbb{E}\left(\frac{\sum_{k=1}^n g_k(M_k)^p}{\sum_{k=1}^n (M_k)^p}\right),$$
$$G^{\vartheta^{(p)}}(n) = \gamma + (1-p) \cdot G_*^{\vartheta^{(p)}}(n).$$

As expected, these formulae reduce to those of (5.10), (5.11) when $p = 1$; and for the Atlas model of (1.7) they give

$$(5.20) \qquad G_*^{\vartheta^{(p)}}(n) = \frac{g}{p} \cdot \left[1 - n \cdot \mathbb{E}\left(\frac{(M_n)^p}{\sum_{k=1}^n (M_k)^p}\right)\right],$$
$$G^{\vartheta^{(p)}}(n) = g + (1-p)G_*^{\vartheta^{(p)}}(n).$$

We shall also look at modified versions of the portfolios considered so far in Examples 5.1–5.3, which "shun the smallest stock in the market."

EXAMPLE 5.4. The restricted market portfolio rule $\widehat{\mu}(\cdot) = (\widehat{\mu}_1(\cdot), \ldots, \widehat{\mu}_n(\cdot))$ is defined in terms of the market portfolio of Example 5.1 as

$$(5.21) \qquad \widehat{\mu}_i(t) := \frac{\mu_i(t)}{1 - \mu_{(n)}(t)} \cdot \mathbb{1}_{\{X_i(t) > X_{(n)}(t)\}}, \qquad i = 1, \ldots, n.$$

We shall justify in the Appendix the computations

$$
\begin{aligned}
(5.22) \qquad G_*^{\widehat{\mu}}(n) &= -\mathbb{E}\left( \frac{g_n M_{n-1} + \sum_{k=1}^{n-1} g_k M_k}{1 - M_n} \right), \\
G^{\widehat{\mu}}(n) &= \gamma - g_n \cdot \mathbb{E}\left( \frac{M_{n-1}}{1 - M_n} \right).
\end{aligned}
$$

These quantities are the same in the case (1.7) of the Atlas model:

$$(5.23) \qquad G_*^{\widehat{\mu}}(n) = G^{\widehat{\mu}}(n) = g\left( 1 - (n-1) \cdot \mathbb{E}\left( \frac{M_{n-1}}{1 - M_n} \right) \right).$$

EXAMPLE 5.5. The restricted equally-weighted portfolio rule $\widehat{\eta}(\cdot) = (\widehat{\eta}_1(\cdot), \ldots, \widehat{\eta}_n(\cdot))$ assigns equal weights to all stocks but the smallest, which receives zero weight:

$$(5.24) \qquad \widehat{\eta}_i(t) := \frac{1}{n-1} \cdot \mathbb{1}_{\{X_i(t) > X_{(n)}(t)\}}, \qquad i = 1, \ldots, n.$$

Clearly,

$$(5.25) \qquad \gamma^{\widehat{\eta}}(t) = \gamma - \frac{g_n}{n-1} + \gamma_*^{\widehat{\eta}}(t), \qquad \gamma_*^{\widehat{\eta}}(t) = \frac{n-2}{2(n-1)^2} \sum_{k=1}^{n-1} \sigma_k^2.$$

For the (generalized) Atlas model of (1.7) with variance structure (4.16), this gives

$$(5.26) \qquad G^{\widehat{\eta}}(n) = G_*^{\widehat{\eta}}(n) = \frac{n(n-2)}{2(n-1)^2}\left( \sigma^2 + s^2 \frac{n-1}{2} \right).$$

Just as in (5.16), we get then

$$(5.27) \qquad \Gamma_*^{\widehat{\eta}} = \widehat{\Gamma}^{\widehat{\eta}} = \frac{\sigma^2}{2}$$

in the case of equal variances ($s = 0$), whereas $\Gamma_*^{\widehat{\eta}} = \widehat{\Gamma}^{\widehat{\eta}} = \infty$ when $s > 0$.

EXAMPLE 5.6. The restricted diversity-weighted portfolio rule $\widehat{\vartheta}^{(p)}(\cdot) = (\widehat{\vartheta}_1^{(p)}(\cdot), \ldots, \widehat{\vartheta}_n^{(p)}(\cdot))$ is defined as

$$
\begin{aligned}
(5.28) \qquad \widehat{\vartheta}_i^{(p)}(t) &:= \frac{(\mu_i(t))^p}{(\mu_{(1)}(t))^p + \cdots + (\mu_{(n-1)}(t))^p} \cdot \mathbb{1}_{\{X_i(t) > X_{(n)}(t)\}}, \\
&\qquad\qquad\qquad\qquad\qquad\qquad\qquad i = 1, \ldots, n,
\end{aligned}
$$

for some fixed number $p \in (0, 1)$. Note that $\widehat{\vartheta}^{(1)}(\cdot)$ is simply the restricted market portfolio $\widehat{\mu}(\cdot)$ of Example 5.4. We shall see in the Appendix that the a.s. limits of (5.6) exist for this portfolio and are given by

$$(5.29) \quad G_*^{\widehat{\vartheta}^{(p)}}(n) = -\frac{1}{p} \cdot \mathbb{E}\left(\frac{g_n(M_{n-1})^p + \sum_{k=1}^{n-1} g_k(M_k)^p}{\sum_{k=1}^{n-1}(M_k)^p}\right),$$

$$(5.30) \quad G^{\widehat{\vartheta}^{(p)}}(n) = \gamma - \frac{1}{p} \cdot \mathbb{E}\left(\frac{g_n(M_{n-1})^p + (1-p)\sum_{k=1}^{n-1} g_k(M_k)^p}{\sum_{k=1}^{n-1}(M_k)^p}\right),$$

in the context of the first-order model. Again, these formulae reduce to those of (5.22) when $p = 1$.

Just as in (5.23), these quantities are the same in the context of (1.7), namely,

$$(5.31) \qquad G_*^{\widehat{\vartheta}^{(p)}}(n) = G^{\widehat{\vartheta}^{(p)}}(n) = \frac{g}{p}\left[1 - (n-1)\cdot\mathbb{E}\left(\frac{(M_{n-1})^p}{\sum_{k=1}^{n-1}(M_k)^p}\right)\right].$$

## 6. Comparisons and approximations.
We can begin now to make some comparisons of long-term-growth behavior for the portfolio rules introduced in Examples 5.1–5.6. For instance, (5.12) and (5.16) give

$$(6.1) \qquad \Gamma^{\mu} = g \quad \text{and} \quad \Gamma_*^{\eta} = (\sigma^2/2), \qquad \Gamma^{\eta} = g + (\sigma^2/2) = g(1 + \alpha)$$

in the context of the Atlas model (1.7), (4.14) and with the notation of (4.15) for the asymptotic (as $n \to \infty$) long-term growth rates of the market portfolio rule $\mu(\cdot)$ and the equally-weighted portfolio rule $\eta(\cdot)$. Thus, *in a large Atlas model, the equally-weighted portfolio outperforms the market portfolio in terms of long-term growth rate, by the amount $\sigma^2/2$.*

Do these features persist when one is not able to invest in the smallest stock? We know from (5.27) that the restricted equally-weighted portfolio $\widehat{\eta}(\cdot)$ has asymptotic long-term growth rate

$$(6.2) \qquad\qquad \Gamma^{\widehat{\eta}} = (\sigma^2/2) = \alpha g = \Gamma_*^{\widehat{\eta}}.$$

*The inability to invest in the smallest stock of a large Atlas model penalizes the long-term growth rate of equal-weighting by the amount $(\sigma^2/2)$, but leaves its long-term excess growth rate the same.*

But how about the performance of the restricted market portfolio $\widehat{\mu}(\cdot)$ of Example 5.4? From (5.23) we see that we have to calculate the limit

$$(6.3) \qquad \Gamma^{\widehat{\mu}} = g\left[1 - \lim_{n\to\infty}\left((n-1)\cdot\mathbb{E}\left(\frac{M_{n-1}}{1-M_n}\right)\right)\right] = \Gamma_*^{\widehat{\mu}},$$

which is not a straightforward task. We work similarly for the diversity-weighted portfolio $\vartheta^{(p)}(\cdot)$ of Example 5.3 and its restricted counterpart $\widehat{\vartheta}^{(p)}(\cdot)$ of Example 5.6: in accordance with (5.20) and (5.31), we have to compute the quantities

$$
\begin{aligned}
(6.4) \qquad \Gamma_*^{\vartheta^{(p)}} &= \frac{g}{p} \cdot \left[ 1 - \lim_{n \to \infty} \left( n \cdot \mathbb{E}\left( \frac{(M_n)^p}{\sum_{k=1}^n (M_k)^p} \right) \right) \right], \\
\Gamma^{\vartheta^{(p)}} &= g + (1-p)\Gamma_*^{\vartheta^{(p)}}
\end{aligned}
$$

and

$$
(6.5) \qquad \Gamma_*^{\widehat{\vartheta}^{(p)}} = \Gamma^{\widehat{\vartheta}^{(p)}} = \frac{g}{p} \left[ 1 - \lim_{n \to \infty} \left( (n-1) \cdot \mathbb{E}\left( \frac{(M_{n-1})^p}{\sum_{k=1}^{n-1} (M_k)^p} \right) \right) \right].
$$

To carry out the computations of (6.3)–(6.5) we shall resort to the certainty-equivalent approximation of Section 4.1. In particular, we shall replace in (6.3)–(6.5) the random variables $M_1, \ldots, M_n$ by the constants

$$
(6.6) \qquad M_k^{CE} \sim \frac{k^{-(\alpha+\beta)} e^{-2\beta k}}{\sum_{j=1}^n j^{-(\alpha+\beta)} e^{-2\beta j}}, \qquad k = 1, \ldots, n \text{ with } \alpha > 0, \beta \geq 0,
$$

as in (4.17) and Examples 4.1, 4.2.

6.1. *Atlas model with* $\alpha = (\sigma^2/2g) > 1$. In this setting we have $\beta = 0$ and the series $\sum_{j \in \mathbb{N}} j^{-\alpha}$ appearing in the denominator of

$$
(6.7) \qquad M_k^{CE} \sim \frac{k^{-\alpha}}{\sum_{j=1}^n j^{-\alpha}}, \qquad k = 1, \ldots, n,
$$

converges, therefore

$$
M_n^{CE} \sim O(n^{-\alpha}), \qquad n M_n^{CE} \sim O(n^{1-\alpha}) \qquad \text{as } n \to \infty.
$$

It follows then from (5.12) that

$$
(6.8) \qquad \Gamma_*^{\mu} \sim g \left( 1 - \lim_{n \to \infty} (n M_n^{CE}) \right) \sim g = \Gamma^{\mu}.
$$

EXAMPLE 5.4 (Continued). For the restricted market portfolio $\widehat{\mu}(\cdot)$ the relations (6.3) and (6.7) now give

$$
(6.9) \qquad \Gamma^{\widehat{\mu}} = \Gamma_*^{\widehat{\mu}} \sim g \left[ 1 - \lim_{n \to \infty} \left( \frac{(n-1) \cdot M_{n-1}^{CE}}{1 - M_n^{CE}} \right) \right] \sim g.
$$

Comparing with (6.1) we see that, in this case, *dropping* (or inability to invest in) *the smallest stock does not result in loss of long-term growth for the market portfolio.*

EXAMPLE 5.3 (Continued). For the diversity-weighted portfolio $\vartheta^{(p)}(\cdot)$ the relations (6.4) and (6.7) give

$$(6.10) \qquad 1 - \frac{p}{g}\Gamma_*^{\vartheta^{(p)}} \sim \lim_{n\to\infty}\left(\frac{n(M_n^{CE})^p}{\sum_{k=1}^n (M_k^{CE})^p}\right) \sim \lim_{n\to\infty}\left(\frac{n^{1-\alpha p}}{\sum_{k=1}^n k^{-\alpha p}}\right).$$

We need to distinguish two cases:

*Case* I. If $(1/\alpha) < p < 1$, then the limit in the expression (6.10) is equal to zero, and we obtain

$$(6.11) \qquad \Gamma_*^{\vartheta^{(p)}} \sim \frac{g}{p}, \qquad \Gamma^{\vartheta^{(p)}} \sim g + (1-p)\Gamma_*^{\vartheta^{(p)}} \sim \frac{g}{p}.$$

Comparing with (6.1) and (6.8) we see an *advantage over the market portfolio* $\mu(\cdot)$ in this case: $\Gamma^{\vartheta^{(p)}} > \Gamma^\mu$, $\Gamma_*^{\vartheta^{(p)}} > \Gamma_*^\mu$. But comparing with (6.1) again, we see that *a disadvantage emerges vis-à-vis the equally-weighted portfolio* $\eta(\cdot)$: $\Gamma^\eta > \Gamma^{\vartheta^{(p)}}$, $\Gamma_*^\eta > \Gamma_*^{\vartheta^{(p)}}$.

*Case* II. If $0 < p \le (1/\alpha)$, then $\lim_{n\to\infty}\left(\frac{n^{1-\alpha p}}{\sum_{k=1}^n k^{-\alpha p}}\right) = 1 - \alpha p$ and (6.10) gives

$$(6.12) \qquad \Gamma_*^{\vartheta^{(p)}} \sim \alpha g = (\sigma^2/2), \qquad \Gamma^{\vartheta^{(p)}} \sim g + (1-p)(\sigma^2/2).$$

There is a definite advantage over the market portfolio ($\Gamma^{\vartheta^{(p)}} > \Gamma^\mu$, $\Gamma_*^{\vartheta^{(p)}} > \Gamma_*^\mu$) and a disadvantage vis-à-vis the equally-weighted portfolio ($\Gamma^\eta > \Gamma^{\vartheta^{(p)}}$, $\Gamma_*^\eta = \Gamma_*^{\vartheta^{(p)}}$).

EXAMPLE 5.6 (Continued). For the restricted diversity-weighted portfolio $\widehat{\vartheta}^{(p)}(\cdot)$ the relations (6.5) and (6.7) give

$$\Gamma_*^{\widehat{\vartheta}^{(p)}} = \Gamma^{\widehat{\vartheta}^{(p)}} \sim \frac{g}{p}\left[1 - \lim_{n\to\infty}\left(\frac{(n-1)^{1-\alpha p}}{\sum_{k=1}^n k^{-\alpha p}}\right)\right].$$

From the preceding discussion we see

$$(6.13) \qquad \Gamma_*^{\widehat{\vartheta}^{(p)}} = \Gamma^{\widehat{\vartheta}^{(p)}} \sim \frac{g}{p} \qquad \text{in case } \alpha p > 1$$

(i.e., no disadvantage at all for dropping the smallest stock), and

$$(6.14) \qquad \Gamma_*^{\widehat{\vartheta}^{(p)}} = \Gamma^{\widehat{\vartheta}^{(p)}} \sim \alpha g = (\sigma^2/2) \qquad \text{in case } \alpha p \le 1$$

[i.e., a drop of $\Gamma^{\vartheta^{(p)}} - \Gamma^{\widehat{\vartheta}^{(p)}} \sim (1-\alpha p)g$ in long-term growth rate, due to inability to invest in the smallest stock].

If $\alpha p > 1$, the restricted diversity-weighted portfolio outperforms the restricted market portfolio [cf. with (6.9)] and underperforms the restricted equally-weighted portfolio [cf. with (6.2)]. When $\alpha p \le 1$, the advantage versus the restricted market portfolio remains, but the disadvantage vis-à-vis the restricted equally-weighted portfolio disappears.

6.2. *Atlas model with $\alpha = (\sigma^2/2g) \le 1$.* In this case we have again $\beta = 0$ in (6.6), so $M_k^{CE} \sim \frac{k^{-\alpha}}{\sum_{j=1}^n j^{-\alpha}}$ satisfies

$$\lim_{n \to \infty} (n M_n^{CE}) \sim \lim_{n \to \infty} \left( \frac{n^{1-\alpha}}{\sum_{k=1}^n k^{-\alpha}} \right) = 1 - \alpha.$$

(a) For the market portfolio $\mu(\cdot)$ we have now

$$(6.15) \qquad \Gamma^\mu = g, \qquad \Gamma_*^\mu \sim g \left( 1 - \lim_{n \to \infty} (n M_n^{CE}) \right) \sim \alpha g = (\sigma^2/2)$$

from (6.1) and (6.8), and for its restricted version $\widehat{\mu}(\cdot)$ equation (6.3) gives

$$(6.16) \qquad \Gamma^{\widehat{\mu}} = \Gamma_*^{\widehat{\mu}} \sim g \cdot \left[ 1 - \lim_{n \to \infty} \left( \frac{(n-1) M_{n-1}^{CE}}{1 - M_n^{CE}} \right) \right] \sim \alpha g = (\sigma^2/2).$$

In other words: when the market portfolio cannot invest in the smallest stock, there is a loss of long-term growth rate $(\Gamma^\mu > \Gamma^{\widehat{\mu}})$ for $\alpha < 1$, whereas $\Gamma^\mu \sim \Gamma^{\widehat{\mu}}$ for $\alpha = 1$.

(b) For the diversity-weighted portfolio $\vartheta^{(p)}(\cdot)$ it is easy to check from (6.4) that the expressions

$$(6.17) \qquad \Gamma_*^{\vartheta^{(p)}} \sim \alpha g = (\sigma^2/2), \qquad \Gamma^{\vartheta^{(p)}} \sim g + (\sigma^2/2)(1-p)$$

of (6.10) and (6.12) prevail again; the same is true of the expressions

$$(6.18) \qquad \Gamma_*^{\widehat{\vartheta}^{(p)}} \sim (\sigma^2/2), \qquad \Gamma^{\widehat{\vartheta}^{(p)}} \sim (\sigma^2/2)$$

of (6.14) for the restricted counterpart $\widehat{\vartheta}^{(p)}(\cdot)$ of $\vartheta^{(p)}(\cdot)$. In this case *the restricted market, equally-weighted and diversity-weighted portfolios $\widehat{\mu}(\cdot)$, $\widehat{\eta}(\cdot)$ and $\widehat{\vartheta}^{(p)}(\cdot)$ have exactly the same long-term-growth performance.*

6.3. *Maximal growth rate in the Atlas model.* All these comparisons beg the obvious question: *What is the maximum long-term growth rate* $\liminf_{T \to \infty} \frac{1}{T} \times Z^\pi(T)$ *from investment, that one can achieve over all possible portfolio rules* $\pi(\cdot)$? In the context of the Atlas model it is natural to guess that the best such rate can be attained by always investing in the smallest, the Atlas stock: namely, that

$$(6.19) \qquad \liminf_{T \to \infty} \frac{1}{T} Z^\pi(T) \le \lim_{T \to \infty} \frac{1}{T} Z^{\pi^*}(T) = ng \qquad \text{a.s.}$$

holds for every portfolio $\pi(\cdot)$, where

$$(6.20) \qquad \pi_i^*(t) := \mathbb{1}_{\{X_i(t) = X_{(n)}(t)\}}, \qquad i = 1, \ldots, n,$$

is the portfolio that invests always and exclusively in the Atlas stock.

This eminently reasonable guess can be justified rigorously when $ng \geq (\sigma^2/2)$ using the theory for portfolio optimization developed in [9], in particular, Example 6.4.2 on page 282 and Theorem 3.10.1 on page 152. It is interesting then to compare the optimal growth rate $G^{\pi^*}(n) = ng$ of (6.19) with the growth rate

$$G^{\eta}(n) = g + \frac{n-1}{n}\frac{\sigma^2}{2}$$

of the equally-weighted portfolio $\eta(\cdot)$ from (5.15); with that of its modified counterpart $\widehat{\eta}(\cdot)$ from (5.26), namely

$$G^{\widehat{\eta}}(n) = g + \frac{n(n-2)}{(n-1)^2}\frac{\sigma^2}{2};$$

and with the asymptotic (as $n \to \infty$) long-term growth rates

$$\Gamma^{\vartheta(p)} = \lim_{n\to\infty} G^{\vartheta(p)}(n) \sim g\left[1 + \frac{1-p}{p}(1 \wedge \alpha p)\right]$$

of the diversity-weighted portfolio, and

$$\Gamma^{\widehat{\vartheta}(p)} = \lim_{n\to\infty} G^{\widehat{\vartheta}(p)}(n) \sim \frac{g}{p}(1 \wedge \alpha p),$$

of its restricted counterpart, from (6.11), (6.12) and (6.13), (6.14), respectively.

The trouble, of course, is that the portfolio of (6.20) is *extremely hard*, if not impossible, to implement in practice—quite in contrast to the portfolios of Examples 5.1–5.6 which can be implemented with relative ease.

6.4. *Generalized Atlas model of* (4.16). Let us consider now the case of variance coefficients of the form (4.16) with $\sigma^2 > 0$. The certainty-equivalent approximation of (6.6) with $\beta = (s^2/4g) > 0$ now has the advantage that the series $\sum_j j^{-(\alpha+\beta)}e^{-2\beta j}$, appearing in the denominator of (6.6), converges for any values of the parameters $\alpha > 0$, $\beta > 0$; this makes the analysis much easier than before.

In particular, it is checked using (6.3)–(6.6) that we have

(6.21)     $\Gamma_*^{\mu} \sim g = \Gamma^{\mu}$   and   $\Gamma^{\widehat{\mu}} \sim g,$     $\Gamma_*^{\widehat{\mu}} \sim g$

for the market portfolio and its modification $\widehat{\mu}(\cdot)$; that

$$\Gamma_*^{\vartheta(p)} \sim (g/p), \qquad \Gamma^{\vartheta(p)} \sim g + (1-p)\Gamma_*^{\vartheta(p)} \sim (g/p)$$

for the diversity-weighted portfolio; and

$$\Gamma_*^{\widehat{\vartheta}(p)} \sim (g/p), \qquad \Gamma^{\widehat{\vartheta}(p)} \sim (g/p)$$

for its modified version. In other words, the diversity-weighted portfolio loses no long-term performance by shunning (or failing to invest in) the smallest stock; and by selecting the parameter $p \in (0, 1)$ sufficiently small, we see that diversity-weighted portfolios can be constructed that have arbitrarily large long-term growth rates—at least if the number of stocks in the market is large.

The assumption of linear growth of variance coefficients with decreasing size captures quite well the actual measurement of stock-price volatilities reported in Figure 5.5 of [2], page 109. This figure plots the smoothed annualized values of $s_k^2$ against rank $k$ in the entire U.S. equity market for the period 1990–1999 (see the discussion in Section 5.1, page 95, of [2], for details of which securities are included). Recalling from (3.10) that $s_k^2 = \sigma_k^2 + \sigma_{k+1}^2$, it is reasonable to make the approximation $\sigma_k^2 \approx s_k^2/2$. Accordingly, Figure 2, which shows the annualized values of $\sigma_k^2$ against $k$ over the same time period, is a scaled version of Figure 5.5 of [2]. In the figure, the variances $\sigma_k^2$ do appear to grow roughly linearly with rank. Using regression, we have estimated the parameters $\sigma^2$ and $s^2$ of (4.16) to be 0.075 and $6.0 \times 10^{-5}$, respectively. Furthermore, we have computed the annualized excess growth rate of the entire U.S. equity market over the same time period (cf. Figure 1 of [3]) to be about 4.4%; motivated by the observation $\Gamma_*^\mu \sim g$ of (6.21) above, we can estimate that $g \sim 0.044$ over this period. Using our estimates for the parameters $\sigma^2$, $s^2$ and $g$, we have plotted the quantities $M_k^{CE}$ of (6.6) in Figure 3, along with the observed capital distribution curve as given in Figure 5.3 of [2]. The two curves are in rough agreement, with the biggest discrepancies occurring for about the 20 highest-ranked stocks.

6.5. *Simulations.* In order to test the performance of the certainty-equivalent approximation, we have used numerical simulations to measure the diversity of Atlas models under a variety of different conditions. In particular, we wish to see if $f^{CE} := f(M_1^{CE}, \ldots, M_n^{CE})$ is a good approximation for the right-hand side $\mathbb{E}[f(M_1, \ldots, M_n)]$ of (4.11), by comparing $f^{CE}$ with the values over time of $T^{-1} \int_0^T f(\mu_{(1)}(t), \ldots, \mu_{(n)}(t)) \, dt$. We have chosen the function $f$ given by $f(x) := \sum_{i=1}^n x_i^p$. This is the $p$th power of the so-called *diversity* function, which is defined in [2] and also as $D(x)$ in Example 5.3 above. We have set $p = 0.5$; results for other values of $p$ in the interval $[0.2, 0.8]$ are similar. In each panel of Figure 4, the horizontal dotted line represents the value of $f^{CE} = f(M_1^{CE}, \ldots, M_n^{CE})$. The three curves in each panel represent the simulated values of the time averages $T^{-1} \int_0^T f(\mu_{(1)}(t), \ldots, \mu_{(n)}(t)) \, dt$ as $T$ varies. The scale on the horizontal axis is in years, assuming 250 iterations per year. Of the three curves in each panel, the solid curve shows the result when the initial values of the weights agree with $(M_1^{CE}, \ldots, M_n^{CE})$. The dashed curve corresponds to initially equal weights, and the dotted curve corresponds to an initial distribution where one weight is very close to 1.0. In the three left-hand side panels, the volatility $\sigma_k^2$ is constant across ranks; the (constant) growth rate $g$ is chosen appropriately in order to model the three cases $\alpha = 0.5$, 1.0 and 1.5. In the corresponding right-hand side panels, the volatility $\sigma_k^2$ now grows linearly with rank $k$, as in (4.16). The parameters $\sigma^2$ and $s^2$, now taken at a daily frequency, have been set as 1/250th of the corresponding annualized values from Section 6.4 above. Once again the parameter $g$ has been selected in order to achieve the desired values of $\alpha$. In each case, the number of stocks $n$ was taken to be 5000, and the simulation was run over 5000 years (1.25 million iterations).

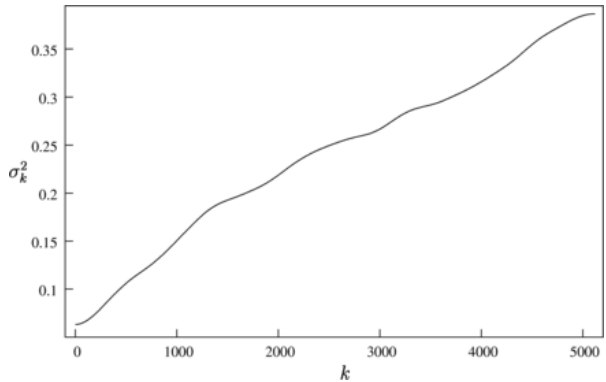

FIG. 2. *Smoothed annualized values of $\sigma_k^2$, from 1990 to 1999 data.*

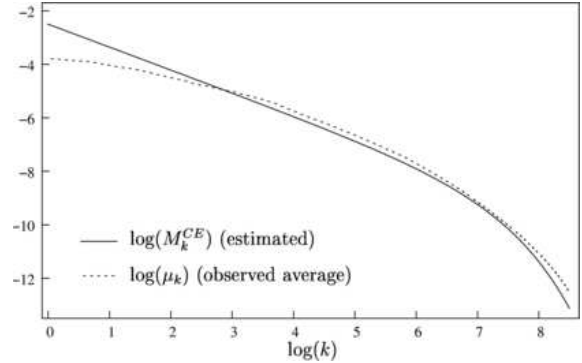

Fig. 3. *Estimated and actual capital distribution curves for* 1990 *to* 1999.

In each simulation, all three initial conditions eventually lead to values of the time-average of $f$ which are relatively stable, approximately equal to each other and which compare favorably with the value arising from the certainty-equivalent approximation. The rate of convergence seems to increase as $\alpha$ decreases. A possible explanation for this is that the largest weight $\mu_{(1)}$ is quite large for higher values of $\alpha$; when this is the case, the diversity $f$ is sensitive to changes in $\mu_{(1)}$, affecting the rate of convergence.

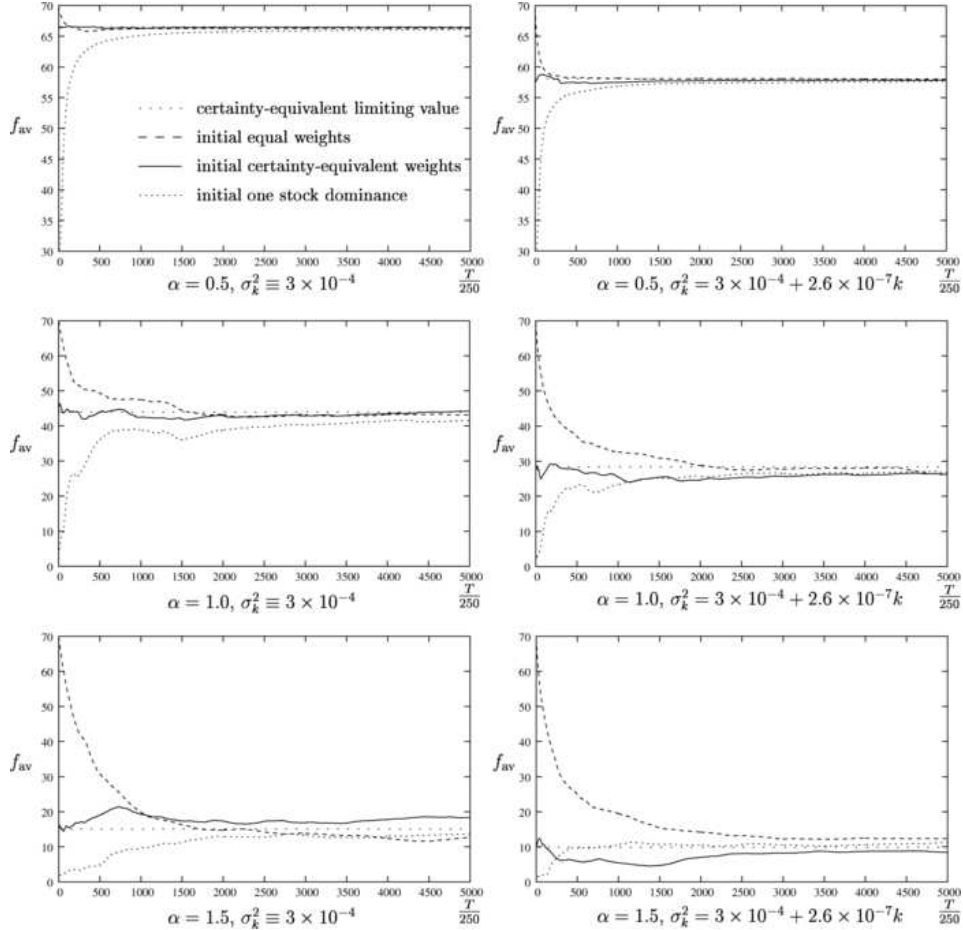

Fig. 4. *Average pth power of diversity in Atlas models ($n = 5000$, $p = 0.5$).*

6.6. *Efficient portfolios in generalized Atlas models.* A portfolio is called *efficient* if its variance is minimal among all portfolios with the same rate of return. Consider a portfolio $\widehat{\pi}(\cdot)$ in a generalized Atlas market with $n + 1$ stocks, restricted to invest in all but the smallest stock at all times; we shall denote by $\widehat{\pi}_{(k)}(t)$ the relative weight this portfolio assigns to the stock ranked $k$th at time $t$. Then from (5.5) we see that the rate of return and the variance of this portfolio are

$$b^{\widehat{\pi}}(t) = \tfrac{1}{2} \sum_{k=1}^{n} \sigma_k^2 \widehat{\pi}_{(k)}(t) \quad \text{and} \quad (\sigma^{\widehat{\pi}}(t))^2 = \sum_{k=1}^{n} \sigma_k^2 (\widehat{\pi}_{(k)}(t))^2,$$

respectively. In order to minimize $(\sigma^{\widehat{\pi}}(t))^2$ with $b^{\widehat{\pi}}(t)$ constant and $\widehat{\pi}_{(1)}(t) \geq 0, \ldots, \widehat{\pi}_{(n)}(t) \geq 0$, $\widehat{\pi}_{(1)}(t) + \cdots + \widehat{\pi}_{(n)}(t) = 1$, we must have

$$2\sigma_k^2 \widehat{\pi}_{(k)}(t) = \lambda_1 \sigma_k^2 + \lambda_2, \qquad k = 1, \ldots, n,$$

where $\lambda_1$ and $\lambda_2$ are Lagrange multipliers. The solution is

$$(6.22) \qquad \widehat{\pi}_{(k)}(t) = \lambda \cdot \frac{1}{n} + (1 - \lambda) \cdot \frac{1}{\sigma_k^2 \sum_{j=1}^{n} \sigma_j^{-2}}, \qquad 0 \leq \lambda \leq 1.$$

The "efficient frontier" for this model consists of the one-parameter family of portfolios defined by (6.22). In the case (4.14) of constant variances the two fractions on the right-hand side of (6.22) are equal, so *for the prototype Atlas model equal weights produce the only efficient portfolio.*

For the general Atlas model, the value $\lambda = 1$ produces the most risky portfolio, and the value $\lambda = 0$ the least risky. With linearly growing variances as in (4.16), *the weights given by* (6.22) *for the portfolios of the efficient frontier are considerably less concentrated in the large stocks than the stable market weights* of Figure 5.3 in [2], page 108.

## 7. Considerations of diversity.

The Atlas model of (1.1), (1.2) has constant and invertible volatility matrix and bounded growth rates, so it admits a unique equivalent martingale measure on *every* finite time-horizon $[0, T]$. For this reason, it cannot be weakly diverse: in other words, for every $T \in (0, \infty)$ and $\delta \in (0, 1)$ we have

$$(7.1) \qquad \mathbb{P}\left( \frac{1}{T} \int_0^T \mu_{(1)}(t)\, dt \leq 1 - \delta \right) < 1,$$

as shown in [4]. We shall argue below that the probability in (7.1) is actually *very close to* 1; for all intents and purposes, such a model captures rather well in practice the descriptive and intuitively plausible requirement, that "no stock should be allowed to dominate the entire market, even on the average, with anything but an extremely low probability."

Let us then try to estimate the probability in (7.1) for a market of $n = 5000$ stocks, similar in size to the U.S. stock market. We take a time period of $T = 2$ years, let $\delta = 0.01$, and wish to measure the probability that the weak diversity condition holds, that is, that none of the stocks has a time-average market weight greater than $1 - \delta = 0.99$ over the 2-year period. We shall assume that the stocks are log-normally distributed relative to the market.

We first need to estimate the relative variance of a typical stock with respect to the market. The drift component, without leakage, of the diversity-weighted portfolio in (5.17) is equal to $(1 - p)/2$, times the weighted average of the relative variances $\tau_{ii}$ of the stocks in the market. In the example of Fernholz [2], Section 6.2, this component was 1.46% per year with $p = 0.50$, so we can estimate the average relative stock variance at about 5.84% per year. This means that the average relative standard deviation would be about 24% per year, which seems within the range of what one would reasonably expect.

For weak diversity to fail, at least one stock must attain a market weight of at least 0.98 at some time during the first year, since otherwise the time-average of the largest weight could not be greater than 0.99 over the 2-year period. By the reflection principle (e.g. [8], Section 2.6.A) the probability of a stock weight hitting 0.98 during the first year is equal to twice the probability that its weight exceeds 0.98 at the end of the year. Let us estimate this last probability.

Suppose the starting weight of the stock is 0.03. For the weight to increase to 0.98, the stock would have to increase relative to the market by a factor $A$, where

$$\frac{0.03A}{0.03A + 0.97} = 0.98,$$

so $A = 1584$. On a logarithmic scale, this would be about $7.37 = \log A$, so with an annual standard deviation of about 0.24, this is slightly greater than 30 standard deviations. If stock price has a log-normal distribution, then the probability of a price move of this size is

$$(7.2) \quad \frac{1}{\sqrt{2\pi}} \int_{30}^{\infty} e^{-t^2/2}\, dt < \frac{1}{30\sqrt{2\pi}} \int_{30}^{\infty} t e^{-t^2/2}\, dt = \frac{e^{-450}}{30\sqrt{2\pi}} \cong 10^{-197}.$$

By the reflection principle, the probability that the stock's market weight attains 0.98 sometime during the year is double this probability. To find the probability that any one of the $n = 5000$ stocks attains a market weight of 0.98, we must multiply the result by 5000, so in all we must multiply the probability in (7.2) by $10^4$. Hence, the probability that at least one of the stocks reaches 0.98 during the first year is not greater than $10^{-193}$, so the probability in (7.1) that the weak diversity condition holds, exceeds $1 - 10^{-193}$. Now even without writing this number out as $0.9999\ldots 999$ we can see that it is pretty close to 1, so it would seem that a market of this type is likely to behave rather like a diverse market over a 2-year period. And this is without invoking antitrust legislation.

## APPENDIX

PROOF OF (2.7). We shall establish in this section the strong law of large numbers (2.7). This property is equivalent to the analogous result

$$(A.1) \qquad \lim_{t \to \infty} \frac{1}{t} Z_k(t) = \gamma \qquad \text{a.s. } \forall\, k = 1, \ldots, n$$

for the log-capitalization processes of (3.1). Indeed, (2.7) implies (A.1) thanks to (3.1) and the partition property (2.2); and conversely, the partition property that precedes (2.2) leads to

$$\sum_{k=1}^{n} \mathbb{1}_{\mathcal{Q}_k^{(j)}}(Y(t)) \cdot Z_k(t) = \sum_{i=1}^{n} Y_i(t) \left( \sum_{k=1}^{n} \mathbb{1}_{\mathcal{Q}_k^{(i)}}(Y(t)) \mathbb{1}_{\mathcal{Q}_k^{(j)}}(Y(t)) \right)$$

$$= \sum_{i=1}^{n} \delta_{ij} Y_i(t) = Y_j(t)$$

for every $j = 1, \ldots, n$, and thus (A.1) leads to (2.7).

We shall prove (A.1) under the assumptions

$$(A.2) \qquad g_1 \leq 0, \qquad g_1 + g_2 \leq 0, \qquad g_1 + \cdots + g_{n-1} \leq 0$$

and

$$(A.3) \qquad\qquad\qquad g_1 + \cdots + g_n = 0$$

which, taken together, are actually weaker than (1.5). To this end, let us recall from (3.3) that

$$(A.4) \quad Z_k(t) = Z_k(0) + (g_k + \gamma)t + \tfrac{1}{2}[\Lambda^{k,k+1}(t) - \Lambda^{k-1,k}(t)] + \sigma_k B_k(t)$$

holds for every $k = 1, \ldots, n$. Fix $k \in \{1, \ldots, n-1\}$ and observe

$$(A.5) \quad \begin{aligned} kZ_k(t) &\leq \sum_{\ell=1}^{k} Z_\ell(t) \\ &= \sum_{\ell=1}^{k} Z_\ell(0) + \left( k\gamma + \sum_{\ell=1}^{k} g_\ell \right) t + \tfrac{1}{2}\Lambda^{k,k+1}(t) + \sum_{\ell=1}^{k} \sigma_\ell B_\ell(t), \end{aligned}$$

$$(A.6) \quad \begin{aligned} (n-k)Z_k(t) &\geq \sum_{\ell=k+1}^{n} Z_\ell(t) \\ &= \sum_{\ell=k+1}^{n} Z_\ell(0) + \left( (n-k)\gamma + \sum_{\ell=k+1}^{n} g_\ell \right) t \\ &\quad - \tfrac{1}{2}\Lambda^{k,k+1}(t) + \sum_{\ell=k+1}^{n} \sigma_\ell B_\ell(t). \end{aligned}$$

After rearranging terms and using (A.3), we see that the two inequalities of (A.5), (A.6) imply

$$(A.7) \quad \begin{aligned} -\frac{1}{2}\left( \frac{1}{k} + \frac{1}{n-k} \right)\Lambda^{k,k+1}(t) &\leq \frac{t}{k}\sum_{\ell=1}^{k} g_\ell - \frac{t}{n-k}\sum_{\ell=k+1}^{n} g_\ell + R_k(t) \\ &= \left( \frac{t}{k} + \frac{t}{n-k} \right)\sum_{\ell=1}^{k} g_\ell + R_k(t), \end{aligned}$$

where

$$(A.8) \quad \begin{aligned} R_k(t) := \frac{1}{k}&\left(\sum_{\ell=1}^{k} Z_\ell(0) + \sum_{\ell=1}^{k} \sigma_\ell B_\ell(t)\right) \\ &- \frac{1}{n-k}\left(\sum_{\ell=k+1}^{n} Z_\ell(0) + \sum_{\ell=k+1}^{n} \sigma_\ell B_\ell(t)\right). \end{aligned}$$

From the law of the iterated logarithm for Brownian motion we observe that we have $\lim_{t\to\infty}(t^{-3/4} R_k(t)) = 0$ a.s., thus

$$(A.9) \quad -\tfrac{1}{2}\Lambda^{k,k+1}(t) \le t\left(\sum_{\ell=1}^{k} g_\ell\right) + o(t^{3/4}) \qquad \text{as } t \to \infty$$

almost surely—for every $k = 1, \ldots, n-1$ thanks to (A.7) and (A.8), and trivially for $k = 0$ and $k = n$.

Let us recall now the decomposition

$$(A.10) \quad Z_k(\cdot) - Z_{k+1}(\cdot) = Z_k(0) - Z_{k+1}(0) + \Theta_k(\cdot) + \Lambda^{k,k+1}(\cdot)$$

of (3.8), which exhibits the nonnegative semimartingale $Z_k(\cdot) - Z_{k+1}(\cdot)$ as the reflection at the origin of the process $\Theta_k(\cdot)$ of (3.9). Thanks to (A.9), the bounded variation part $(g_k - g_{k+1})t - \frac{1}{2}[\Lambda^{k-1,k}(t) + \Lambda^{k+1,k+2}(t)]$ of this semimartingale $\Theta_k(\cdot)$ is dominated by

$$(g_k - g_{k+1})t + \left(\sum_{\ell=1}^{k-1} g_\ell + \sum_{\ell=1}^{k+1} g_\ell\right)t + o(t^{3/4}) = 2t\left(\sum_{\ell=1}^{k} g_\ell\right) + o(t^{3/4}) \qquad \text{as } t \to \infty$$

a.s. In other words, $\Theta_k(\cdot)$ is dominated by a Brownian motion [thus $Z_k(\cdot) - Z_{k+1}(\cdot)$ is dominated by a reflected Brownian motion] whose drift is bounded from above by $o(t^{-1/4})$ a.s. [If the inequalities in (A.2) are strict, as they are in (1.5), then this drift is eventually negative and we can appeal to the strong law of large numbers, as opposed to the law of the iterated logarithm, for the Brownian motion process.] A reflected Brownian motion with drift of this type is easily seen to be of the order $o(t)$ as $t \to \infty$, and so we have

$$(A.11) \quad \lim_{t\to\infty} \frac{1}{t}(Z_k(t) - Z_{k+1}(t)) = 0 \qquad \text{a.s.}$$

Now let us divide by $t$ throughout (A.10), then let $t \to \infty$ to obtain

$$\lim_{t\to\infty} \frac{1}{t}\left[\frac{1}{2}(\Lambda^{k-1,k}(t) + \Lambda^{k+1,k+2}(t)) - \Lambda^{k,k+1}(t)\right] = g_k - g_{k+1} \qquad \text{a.s.}$$

for every $k = 1, \ldots, n-1$, from the strong law of large numbers for Brownian motion in conjunction with (A.11). Adding up from $k = \ell$ to $k = n - 1$ we get

$$(A.12) \quad \lim_{t\to\infty} \frac{1}{t}[\Lambda^{\ell-1,\ell}(t) - \Lambda^{\ell,\ell+1}(t) - \Lambda^{n-1,n}(t)] = 2(g_\ell - g_n) \qquad \text{a.s.}$$

for any $\ell = 1, \ldots, n-1$, since $\Lambda^{n,n+1}(\cdot) \equiv 0$. Adding up over these values of $\ell$ we obtain

$$\lim_{t \to \infty} \frac{1}{t}[\Lambda^{0,1}(t) - \Lambda^{n-1,n}(t) - (n-1) \cdot \Lambda^{n-1,n}(t)]$$

$$= 2\sum_{\ell=1}^{n-1} g_\ell - 2(n-1) \cdot g_n \qquad \text{a.s.},$$

and recalling (A.3) and $\Lambda^{0,1}(\cdot) \equiv 0$ we arrive at

$$\lim_{t \to \infty} \frac{1}{t}\Lambda^{n-1,n}(t) = 2g_n \qquad \text{a.s.}$$

Substitution into (A.12) yields

$$\lim_{t \to \infty} \frac{1}{t}[\Lambda^{\ell-1,\ell}(t) - \Lambda^{\ell,\ell+1}(t)] = 2g_\ell \qquad \text{a.s.}$$

for $\ell = 1, \ldots, n$. The property (A.1) now follows from this, in conjunction with (A.4) and the strong law of large numbers for Brownian motion. $\square$

REMARK. Suppose now that the condition (A.2) fails; namely, that $\sum_{\ell=1}^k g_\ell > 0$ holds for some $k = 1, \ldots, n-1$. Then by analogy with (A.5), (A.6) we have

$$kZ_1(t) \geq \sum_{\ell=1}^k Z_\ell(t) = \sum_{\ell=1}^k Z_\ell(0) + \left(k\gamma + \sum_{\ell=1}^k g_\ell\right)t + \tfrac{1}{2}\Lambda^{k,k+1}(t) + \sum_{\ell=1}^k \sigma_\ell B_\ell(t),$$

and the strong law of large numbers for Brownian motion implies

$$\liminf_{t \to \infty} \frac{Z_1(t)}{t} \geq \gamma + \frac{1}{k}\sum_{\ell=1}^k g_\ell > \gamma \qquad \text{a.s.}$$

so that (A.1) cannot hold in this case.

PROOF OF (2.18). To see that the process $\widetilde{Y}(\cdot)$ is recurrent with respect to $B \cap \Pi$ for some ball $B \subset \mathbb{R}^n$ centered at the origin, it suffices to show that the process $N(\cdot) := (1 + \|\widetilde{Y}(\cdot)\|^2)^{1/2}$ is recurrent with respect to $[1, C]$ for some $C > 0$. An application of Itô's rule to (2.15) shows that

(A.13)
$$\begin{aligned} dN(t) = {} & ((N(t))^{-1}\widetilde{Y}(t) \cdot G(\widetilde{Y}(t)) + N_1(t)) \, dt \\ & + ((N(t))^{-1}\widetilde{Y}(t) \cdot \widetilde{S}(\widetilde{Y}(t)) \, dW(t)), \end{aligned}$$

where

$$N_1(t) := (N(t))^{-1}\left(1 - \frac{1}{n}\right)\sum_{i=1}^n \sigma_i^2 - (N(t))^{-3}\sum_{i=1}^n (\widetilde{Y}_{p_t(i)}(t))^2 \sigma_i^2.$$

Here $p_t$ is defined as in (1.4); that is, $p_t \in \Sigma_n$ satisfies $\widetilde{Y}(t) \in \mathcal{R}_{p_t}$ for all $t \geq 0$. It is easy to check that $N_1(t) \leq c_1(N(t))^{-1}$ for some $c_1 > 0$; furthermore, the fact that $\widetilde{Y}(t) \cdot G(\widetilde{Y}(t)) \leq c \|\widetilde{Y}(t)\|$ for some $c < 0$ implies that the term $(N(t))^{-1}\widetilde{Y}(t) \cdot G(\widetilde{Y}(t))$ in (A.13) is uniformly bounded from above by some constant $c_2 < 0$ whenever $N(t) \geq 2$. It follows that there are constants $c_3 < 0$, $C > 0$ such that the drift term in (A.13) is bounded from above by $c_3$ whenever $N(t) > C$. Finally, note that the coefficient $(N(t))^{-1}\widetilde{Y}(t) \cdot \widetilde{S}(\widetilde{Y}(t))$ of $dW(t)$ in (A.13) is a matrix whose entries are uniformly bounded from above and below by $\max_k\{\sigma_k\}$ and $\min_k\{\sigma_k\}$, respectively. The desired recurrence follows. $\square$

PROOF OF (2.21). In order to prove that the processes $\widetilde{Y}(\cdot)$ and $p\widetilde{Y}(\cdot)$ of (2.15) and (2.20) respectively [with $\tilde{y} = 0$ in (2.15)] have the same stationary distribution $\mu$, it suffices to show that the processes have the same law. We claim that $\widetilde{Y}(\cdot)$ also satisfies (2.20), that is,

$$(A.14) \qquad d\widetilde{Y}(t) = G_1(\widetilde{Y}(t))\,dt + \widetilde{S}_1(\widetilde{Y}(t))\,dW(t), \qquad \widetilde{Y}(0) = 0.$$

Indeed, if $\widetilde{Y}(\cdot)$ and $p\widetilde{Y}(\cdot)$ both satisfy (A.14), then they have the same law, since the stochastic differential equation (A.14) has a unique solution in the sense of the probability law. To establish (A.14) for the process $\widetilde{Y}(\cdot)$ of (2.15), note that $G(y) = G_1(y)$ and $S(y) = S_1(y)$ except on the set $\mathcal{R}^e = \{(y_1, \ldots, y_n) \in \mathbb{R}^n \mid y_i - y_j = 0 \text{ for some } i \neq j\}$. It suffices to show that

$$(A.15) \qquad \operatorname{meas}\{t \in [0, \infty) : \widetilde{Y}_i(t) - \widetilde{Y}_j(t) = 0\} = 0 \qquad \text{a.s.}$$

for all pairs $i \neq j$. For such pairs $(i, j)$, the process $Y_{ij}(\cdot) := \widetilde{Y}_i(\cdot) - \widetilde{Y}_j(\cdot)$ is a semimartingale, with drift bounded in absolute value by $\max_k\{2|g_k|\}$ and variance bounded from above and below by $\max_k\{2\sigma_k^2\}$ and $\min_k\{2\sigma_k^2\} > 0$, respectively. Over any bounded time-interval $[0, T]$, removal of drift via the Girsanov theorem, followed by a time change, establishes (A.15) and completes the proof of (2.21). $\square$

PROOF OF (5.30), (5.29), (5.22). To justify the computations (5.30), (5.29) we recall from (4.3.4) of [2] that the a.s. identity

$$\log \frac{Z^{\widehat{\vartheta}^{(p)}}(T)}{Z^\mu(T)} = \frac{1}{p} \log \left( \frac{\sum_{k=1}^{n-1}(\mu_{(k)}(T))^p}{\sum_{k=1}^{n-1}(\mu_{(k)}(0))^p} \right)$$

$$+ (1-p) \int_0^T \gamma_*^{\widehat{\vartheta}^{(p)}}(t)\,dt - \frac{1}{2}\int_0^T \widehat{\vartheta}_{(n-1)}^{(p)}(t)\,d\Lambda^{n-1,n}(t)$$

holds. From (5.3) we have

$$\gamma^{\widehat{\vartheta}^{(p)}}(t) - \gamma_*^{\widehat{\vartheta}^{(p)}}(t) = \gamma + \sum_{k=1}^{n-1} g_k \widehat{\vartheta}_{(k)}^{(p)}(t) = \gamma + \frac{\sum_{k=1}^{n-1} g_k(\mu_{(k)}(t))^p}{\sum_{k=1}^{n-1}(\mu_{(k)}(t))^p},$$

which allows us to write the above identity in the equivalent form

$$
\begin{aligned}
\frac{1}{T}&\left(\log Z^{\widehat{\vartheta}(p)}(T) - \int_0^T \gamma^{\widehat{\vartheta}(p)}(t)\,dt\right) \\
&= \frac{1}{T}(\log Z^\mu(T) - \gamma) - \frac{p}{T}\int_0^T \gamma_*^{\widehat{\vartheta}(p)}(t)\,dt \\
&\quad + \frac{1}{pT}\log\left(\frac{\sum_{k=1}^{n-1}(\mu_{(k)}(T))^p}{\sum_{k=1}^{n-1}(\mu_{(k)}(0))^p}\right) \\
&\quad - \frac{1}{T}\int_0^T \frac{\sum_{k=1}^{n-1} g_k(\mu_{(k)}(t))^p}{\sum_{k=1}^{n-1}(\mu_{(k)}(t))^p}\,dt - \frac{1}{2T}\int_0^T \widehat{\vartheta}_{(n-1)}^{(p)}(t)\,d\Lambda^{n-1,n}(t).
\end{aligned}
$$
(A.16)

On the other hand, we get from (3.5), (3.7), (4.11) and (1.5) that

$$
\lim_{T\to\infty}\frac{1}{2T}\int_0^T \widehat{\vartheta}_{(n-1)}^{(p)}(t)\,d\Lambda^{n-1,n}(t) = g_n \cdot \mathbb{E}\left(\frac{(M_{n-1})^p}{\sum_{k=1}^{n-1}(M_k)^p}\right)
$$
(A.17)

holds a.s. Taking now the limit as $T$ goes to infinity in (A.16) and using (5.4), (5.10), (4.11) and (A.17), we see that the a.s. limits in (5.6) exist for this portfolio and are given by the formulae of (5.30), (5.29). Finally, (5.22) is simply a special case of (5.30), (5.29) in the case $p=1$ of the market portfolio. □

**Note added in proof.** We show that the long-term average relative capitalization weight, for each individual stock (listed by name) in a first-order model, tends to $1/n$: in other words, for each $i=1,\ldots,n$, we have

$$
\lim_{T\to\infty}\frac{1}{T}\int_0^T \mu_i(t)\,dt = \frac{1}{n} \qquad \text{a.s.}
$$
(N.1)

It is not clear that the left-hand limit exists. To show that it does, first set

$$
g(y_1,\ldots,y_n) := \left(\frac{\exp(y_1)}{\sum_{i=1}^n \exp(y_i)},\ldots,\frac{\exp(y_n)}{\sum_{i=1}^n \exp(y_i)}\right), \qquad (y_1,\ldots,y_n)\in\mathbb{R}^n.
$$

By (2.19), we have for any $i,k=1,\ldots,n$ and bounded, measurable $f$,

$$
\begin{aligned}
\lim_{T\to\infty}\frac{1}{T}&\int_0^T f\circ g(\widetilde{Y}(t))\mathbb{1}_{\mathcal{Q}_k^{(i)}}(\widetilde{Y}(t))\,dt \\
&= \int_\Pi f\circ g(y)\mathbb{1}_{\mathcal{Q}_k^{(i)}}(y)\,d\mu(y) \qquad \text{a.s.}
\end{aligned}
$$
(N.2)

Suppose that $f$ is symmetric in all variables and let $p\in\Sigma_n$ be a permutation of $\{1,\ldots,n\}$ such that $p(j)=i$ for some $j\in\{1,\ldots,n\}$. Since the measure $\mu$ is invariant under the action of $\Sigma_n$, $f\circ g(py)=f\circ g(y)$, and $\mathbb{1}_{\mathcal{Q}_k^{(i)}}(py)=\mathbb{1}_{\mathcal{Q}_k^{(i)}}(y)$ for $y$ not in the $\mu$-null set $\mathcal{R}^e$ [defined after (2.11)], we have

$$
\int_\Pi f\circ g(y)\mathbb{1}_{\mathcal{Q}_k^{(i)}}(y)\,d\mu(y) = \int_\Pi f\circ g(y)\mathbb{1}_{\mathcal{Q}_k^{(j)}}(y)\,d\mu(y).
$$

It follows that the almost sure limit on the left-hand side of (N.2) is independent of $i$. By the definition of $g$ and the symmetry of $f$, this limit may be expressed as

$$\lim_{T \to \infty} \frac{1}{T} \int_0^T f(\mu_{(1)}(t), \ldots, \mu_{(n)}(t)) \mathbb{1}_{\mathcal{Q}_k^{(i)}}(Y(t)) \, dt.$$

Since $\{\mathcal{Q}_k^{(i)}\}_{1 \le i \le n}$ is a partition of $\mathbb{R}^n$ for each fixed $k$, summation of the previous expression over $i$ and (4.11) lead to the following refinement of (4.11):

$$(4.11)' \qquad \begin{aligned} \lim_{T \to \infty} \frac{1}{T} &\int_0^T f(\mu_{(1)}(t), \ldots, \mu_{(n)}(t)) \mathbb{1}_{\mathcal{Q}_k^{(i)}}(Y(t)) \, dt \\ &= \frac{1}{n} \mathbb{E}[f(M_1, \ldots, M_n)] \qquad \text{a.s.} \end{aligned}$$

In particular, with $f(y_1, \ldots, y_n) = y_k$ on the set $\{y \in \mathbb{R}^n | y_1 \ge \cdots \ge y_n \ge 0, \sum_{i=1}^n y_i = 1\}$,

$$\lim_{T \to \infty} \frac{1}{T} \int_0^T \mu_i(t) \mathbb{1}_{\mathcal{Q}_k^{(i)}}(Y(t)) \, dt = \lim_{T \to \infty} \frac{1}{T} \int_0^T \mu_{(k)}(t) \mathbb{1}_{\mathcal{Q}_k^{(i)}}(Y(t)) \, dt$$

$$= \frac{1}{n} \mathbb{E}[M_k] \qquad \text{a.s.}$$

Summation over $k$, (2.2), and the fact that $\sum_{k=1}^n M_k \equiv 1$ lead to the desired result (N.1).

**Acknowledgment.** The authors are indebted to Professor Ruth Williams for her very helpful comments and corrections on an earlier version of this paper.

A. D. BANNER
R. FERNHOLZ
INTECH
ONE PALMER SQUARE
PRINCETON, NEW JERSEY 08542
USA
E-MAIL: adrian@enhanced.com
E-MAIL: bob@enhanced.com

I. KARATZAS
DEPARTMENTS OF MATHEMATICS
   AND STATISTICS
COLUMBIA UNIVERSITY
NEW YORK, NEW YORK 10027
USA
E-MAIL: ik@math.columbia.edu